\documentclass[12pt]{amsart}
\usepackage{amsmath}
\usepackage{amssymb}
\usepackage{amsfonts}
\usepackage{amsthm}
\usepackage{amscd}
\usepackage{bm}
\newtheorem{thrm}{Theorem}

\newtheorem{prop}{Proposition}

\newtheorem{thm}{Theorem}[section]
\newtheorem{lem}[thm]{Lemma}
\newtheorem{pro}[thm]{Proposition}
\newtheorem{cor}[thm]{Corollary}

\newcommand{\C}{\mathbf{C}}

\newcommand{\Z}{\mathbf{Z}}
\newcommand{\Q}{\mathbf{Q}}
\newcommand{\N}{\mathbf{N}}

\begin{document}
\title[On the graded quotients of Fricke characters]
{On the graded quotients of \\ the ring of Fricke characters of a free group}
\address[Eri Hatakenaka]{Tokyo University of Agriculture and Technology,
        2-24-16, Naka-cho, Koganei-shi, Tokyo, 184-8588, Japan}
\email{hataken@cc.tuat.ac.jp}
\address[Takao Satoh]{Department of Mathematics, Faculty of Science Division II, Tokyo University of Science,
         1-3 Kagurazaka, Shinjuku, Tokyo, 162-8601, Japan}
\email{takao@rs.tus.ac.jp}
\maketitle
\begin{center}
{\textsf{Dedicated to Professor Sadayoshi Kojima on the occasion of his 60th birthday}}
\end{center}

\vspace{0.5em}

\begin{center}
{\sc Eri Hatakenaka}\footnote{e-address: hataken@cc.tuat.ac.jp} and {\sc Takao Satoh}\footnote{e-address: takao@rs.tus.ac.jp} \\
\end{center}

\begin{abstract}
In this paper, for a group $G$, we consider the $\mathrm{Aut}\,G$-invariant ideal $J$ generated by $\mathrm{tr}\,x -2$ for all $x \in G$ in the 
ring of Fricke characters of $G$. We study a descending filtration $J \supset J^2 \supset J^3 \supset \cdots$, and
its graded quotients $\mathrm{gr}^k(J):=J^k/J^{k+1}$ for $k \geq 1$.
The first purpose of the paper is to determine the $\Q$-vector space structures of $\mathrm{gr}^k(J)$ in the case where $G$ is a free group
$F_n$ of rank $n \geq 2$, and $k=1$ and $2$.

Next, we introduce a normal subgroup $\mathcal{E}_G(k)$ consisting of automorphisms of $G$ which act on $J/J^{k+1}$ trivially.
These normal subgroups define a central filtration of $\mathrm{Aut}\,G$. This is a Fricke character analogue of the Andreadakis-Johnson filtration
$\mathcal{A}_G(k)$ of $\mathrm{Aut}\,G$.
The main purpose of the paper is to show that $\mathcal{E}_{F_n}(1)$ is equal to $\mathrm{Inn}\,F_n \cdot \mathcal{A}_{F_n}(2)$
where $\mathrm{Inn}\,F_n$ is the inner automorphism group of a free group $F_n$, and
that $\mathcal{A}_{F_n}(2k) \subset \mathcal{E}_{F_n}(k)$ for any $k \geq 1$.
\end{abstract}
\section{Introduction}\label{S-Int}

Let $G$ be a group generated by elements $x_1, \ldots, x_n$. We denote by
\[ R(G) := \mathrm{Hom}(G, \mathrm{SL}(2,\C)) \]
the set of all group homomorphisms from $G$ to $\mathrm{SL}(2,\C)$. Let
\[ \mathcal{F}(R(G),\C) := \{ \chi : R(G) \rightarrow \C \} \]
be the set of all complex-valued functions of $R(G)$. Then we can consider $\mathcal{F}(R(G),\C)$ as a commutative ring in a natural way.
For any $x \in G$, we define an element $\mathrm{tr}\,x \in \mathcal{F}(R(G),\C)$ to be
\[ (\mathrm{tr}\,x)(\rho) := \mathrm{tr}\, \rho(x) \]
for any $\rho \in R(G)$. Here \lq\lq $\mathrm{tr}$" in the right hand side means the trace of $2 \times 2$ matrix $\rho(x) \in \mathrm{SL}(2,\C)$.
The element $\mathrm{tr}\,x$ in
$\mathcal{F}(R(G),\C)$ is called the Fricke character of $x \in G$. 
Let $\mathfrak{X}(G)$ be the $\Z$-submodule of $\mathcal{F}(R(G),\C)$ generated by all $\mathrm{tr}\,x$ for $x \in G$. Then 
$\mathfrak{X}(G)$ is closed under the multiplication of $\mathcal{F}(R(G),\C)$. (See Subsection {\rmfamily \ref{Ss-Des}} for details.)

\vspace{0.5em}

Classically, Fricke characters were begun to studied by Fricke for a free group $F_n$ on $x_1, \ldots, x_n$ in connection with certain problems
in the theory of Riemann surfaces. (See \cite{Fri}.) In 1970, Horowitz \cite{Ho1} and \cite{Ho2} investigated algebraic properties of
$\mathfrak{X}(F_n)$ using the combinatorial group theory. In particular, he \cite{Ho1} showed that for any $x \in F_n$, the Fricke character $\mathrm{tr}\,x$
can be written as a polynomial with integral coefficients in $2^n-1$ characters
$\mathrm{tr}\, x_{i_1} x_{i_2} \cdots x_{i_l}$ for $1 \leq l \leq n$ and $1 \leq i_1 < i_2 < \cdots < i_l \leq n$. 
He \cite{Ho2} also showed that the subgroup of $\mathrm{Aut}\,F_n$
consisting of automorphisms which act on $\mathfrak{X}(F_n)$ tirivially is just the inner automorphism group $\mathrm{Inn}\,F_n$ of $F_n$. Namely,
the action of $\mathrm{Aut}\,F_n$ on the ring of Fricke characters $\mathfrak{X}(F_n)$ induces a faithful
representatrion of the outer automorphism group $\mathrm{Out}\,F_n := \mathrm{Aut}\,F_n/\mathrm{Inn}\,F_n$. 
However, since $\mathfrak{X}(F_n)$ is not finitely generated as a $\Z$-module, it is difficult to handle this representation
directly.

\vspace{0.5em}

On the other hand, in order to make the structure of the Fricke characters $\mathfrak{X}(F_n)$ clear,
it is important to study the ideal of polynomials
in the characters which vanishes on any representations of $G$. More precisely, consider a polynomial ring
\[ \Z[t] := \Z[t_{i_1 \cdots i_l} \,|\, 1 \leq l \leq n, \,\, 1 \leq i_1 < i_2 < \cdots < i_l \leq n] \]
of $2^n-1$ indeterminates, and an ideal
\[ I = \{ f \in \Z[t] \,|\, f(\mathrm{tr}\,\rho(x_{i_1} \cdots x_{i_l}))=0 \,\,\, \text{for any} \,\,\, \rho \in R(G) \}. \]
In \cite{Ho1}, for $G=F_n$, Horowitz showed that $I$ is trivial for $n=1$ and $2$, and is principal for $n=3$.
Whittemore \cite{Whi} showed that $I$ is not principal for $G=F_n$ and $n \geq 4$. Although the ideal $I$ has been studied by many authors
for over forty years, very little is known for it.

\vspace{0.5em}

Here, we consider the rationalization of the situation above.
Let $\mathfrak{X}_{\Q}(G)$ be the $\Q$-subspace of $\mathcal{F}(R(G),\C)$ generated by $\mathrm{tr}\,x$ for any $x \in G$.
Similary to $\mathfrak{X}(G)$, $\mathfrak{X}_{\Q}(G)$ is closed under the multiplication of $\mathcal{F}(R(G),\C)$, and has a multiplicative
unit $1 = \frac{1}{2} \mathrm{tr}\, 1_G$. Hence, $\mathfrak{X}_{\Q}(G)$ is a ring.
We call $\mathfrak{X}_{\Q}(G)$ the ring of Fricke characters of $G$ over $\Q$.
By observing the formula (\ref{eq-5}), which is given in Subsection {\rmfamily \ref{Ss-Des}} below, and Horowitz's result, we see that
for any $x \in G$, the Fricke character $\mathrm{tr}\,x$
can be written as a polynomial with ratinal coefficients in $n + \binom{n}{2} + \binom{n}{3}$ characters
$\mathrm{tr}\, x_{i_1} x_{i_2} \cdots x_{i_l}$ for $1 \leq l \leq 3$ and $1 \leq i_1 < i_2 < \cdots < i_l \leq n$. 
Consider a polynomial ring
\[ \Q[t] := \Q[t_{i_1 \cdots i_l} \,|\, 1 \leq l \leq 3, \,\, 1 \leq i_1 < i_2 < \cdots < i_l \leq n] \]
and its ideal
\[ I_{\Q} := \{ f \in \Q[t] \,|\, f(\mathrm{tr}\,\rho(x_{i_1} \cdots x_{i_l}))=0 \,\,\, \text{for any} \,\,\, \rho \in R(G) \}. \]
Similarly to $I$, the ideal $I_{\Q}$ plays important roles in the various study of the ring structure of $\mathfrak{X}_{\Q}(G)$.
One of the most advantages to consider the rationalization of the Fricke characters is that the number of the indeterminates of $\Q[t]$
is fewer than that of $\Z[t]$, and it makes various computation much easy to handle.

\vspace{0.5em}

In the present paper, in order to construct finite dimensinal representations of $\mathrm{Aut}\,G$,
we consider a descending filtration of 
$\mathrm{Aut}\,G$-invariant ideals of $\Q[t]/I_{\Q}$, and take its graded quotients.
Set $t_{i_1 \cdots i_l}' := t_{i_1 \cdots i_l} -2 \in \Q[t]$. We also denote by $t_{i_1 \cdots i_l}'$ its coset class in $\Q[t]/I_{\Q}$.
Consider an ideal
\[ J := (t_{i_1 \cdots i_l}' \,|\, 1 \leq l \leq 3, \,\, 1 \leq i_1 < i_2 < \cdots < i_l \leq n) \subset \Q[t]/I_{\Q} \]
generated by all $t_{i_1 \cdots i_l}'$'s. Then, we have a descending filtration
\[ J \supset J^2 \supset J^3 \supset \cdots \]
of $\mathrm{Aut}\,G$-invariant ideals of $\Q[t]/I_{\Q}$. (See Subsection {\rmfamily \ref{Ss-Des}} for details.)
Set
\[ \mathrm{gr}^k(J) := J^k/J^{k+1}. \]
Each of $\mathrm{gr}^k(J)$ is $\mathrm{Aut}\,G$-invariant $\Q$-vector space of finite dimension for any $k \geq 1$.
This technique is deeply inspired by a result of Magnus \cite{Mg2} who originally studied the behavior of the action of $\mathrm{Aut}\,F_3$ on
$\mathrm{gr}^1(J)$. In \cite{Mg2}, he pointed out the difficulties to find $\mathrm{Aut}\,F_n$-invariant ideals of $\mathfrak{X}(F_n)$ and
its quotient rings as a finite dimensional representation of $\mathrm{Aut}\,F_n$ in general. Moreover, he \cite{Mg2} also stated that
in order to get accessible situation, it seems to be better to use rational functions rather than integral polynomials.
In this paper, however, we consider the rational polynomials to obtain finite dimensional representations of $\mathrm{Aut}\,F_n$.

\vspace{0.5em}

The first purpose of the paper is to determine the structure of $\mathrm{gr}^k(J)$ for $G=F_n$, $n \geq 2$ and $k=1, 2$. 
Set
\[ T := \{ t_i' \,|\, 1 \leq i \leq n \} \cup \{ t_{ij}' \,|\, 1 \leq i < j \leq n \} \cup \{ t_{ijk}' \,|\, 1 \leq i < j < k \leq n \} \subset J \]
and
\[\begin{split}
   S := & \{ t_i' t_j' \,|\, 1 \leq i \leq j \leq n \} \cup \{ t_i' t_{ab}' \,|\, 1 \leq i \leq n, \,\, 1 \leq a < b \leq n \} \\
        & \cup \{ t_i' t_{abc}' \,|\, 1 \leq i \leq n, \,\, 1 \leq a<b<c \leq n \} \\
        & \cup \{ t_{ij}' t_{ab}' \,|\, 1 \leq i<j \leq n, \,\, 1 \leq a<b \leq n, \,\, (i,j) \leq (a,b) \}, \\
        & \cup \{ t_{ab}' t_{abc}', \, t_{ac}' t_{abc}', \, t_{bc}' t_{abc}' \,|\, 1 \leq a<b<c \leq n \} \\
        & \cup \{ t_{ia}' t_{abc}', t_{ib}' t_{abc}', t_{ic}' t_{abc}', t_{ia}' t_{ibc}', t_{ab}' t_{iac}', t_{ab}' t_{ibc}', t_{ac}' t_{ibc}', t_{ib}' t_{iac}' \,|\, 1 \leq i < a< b<c \leq n \} \\
        & \cup \{ t_{ja}' t_{ibc}', t_{jb}' t_{iac}', t_{jc}' t_{iab}', t_{ab}' t_{ijc}', t_{ac}' t_{ijb}', t_{bc}' t_{ija}'
            \,|\, 1 \leq i<j<a<b<c \leq n \} \\
      \subset & \,\,J^2
  \end{split}\]
respectively. We show
\begin{thrm}($=$ Propositions {\rmfamily \ref{P-1}}, {\rmfamily \ref{P-haruhi}} and {\rmfamily \ref{P-2}}.)
For $G=F_n$ and $n \geq 2$, the sets $T$ and $S$ are basis of the $\Q$-vector spaces $\mathrm{gr}^1(J)$ and $\mathrm{gr}^2(J)$ respectively.
\end{thrm}
In general, it seems to be very complicated to find a basis of $\mathrm{gr}^k(J)$ for general $k \geq 3$.

\vspace{1em}

Next, for any group $G$, we consider a descending filtration of $\mathrm{Aut}\,G$. 
For any $k \geq 1$, let $\mathcal{E}_{G}(k)$ be the subgroup of $\mathrm{Aut}\,G$ consisting of automorphisms which act on
$J/J^{k+1}$ trivially. Then we see that the groups $\mathcal{E}_G(k)$ define a descending filtration
\[ \mathcal{E}_G(1) \supset \mathcal{E}_G(2) \supset \cdots \supset \mathcal{E}_G(k) \supset \cdots \]
of $\mathrm{Aut}\,G$. 

\vspace{0.5em}

This filtration is a Fricke character analogue of the Andreadakis-Johnson filtration $\mathcal{A}_G(k)$ of $\mathrm{Aut}\,G$.
The Andreadakis-Johnson filtration was originally introduced by Andreadakis \cite{And} in 1960's.
The name \lq\lq Johnson" comes from Dennis Johnson who studied this type of filtration for the
mapping class group of a surface in 1980's. It is called the Johnson filtration of the
mapping class group.
The Johnson homomorphisms are originally introduced by Johnson
in order to investigate the graded quotients of the Johnson filtration
in a series of his pioneering works \cite{Jo1}, \cite{Jo2}, \cite{Jo3} and \cite{Jo4}.
In \cite{Mo1}, Morita began to study the Johnson homomorphisms of the mapping class groups and $\mathrm{Aut}\,F_n$
systematically.
Today, together with the theory of the Johnson
homomorphisms, the Andreadskis-Johnson filtration is one of powerful tools to study the group structure of the automorphism group of a group.
(See Section {\rmfamily \ref{S-John}} for notation, and see \cite{S06} or \cite{S14} for basic materials concerning
the Andreadakis-Johnson filtration and the Johnson homomorphisms of $\mathrm{Aut}\,F_n$.)

\vspace{0.5em}

The main purpose of the paper is to show
\begin{prop}\label{I-P-2}($=$ Proposition {\rmfamily \ref{P-3}}.)
For any $k, l \geq 1$, $[\mathcal{E}_G(k), \mathcal{E}_G(l)] \subset \mathcal{E}_G(k+l)$.
\end{prop}
and
\begin{thrm}\label{I-T-1}($=$ Theorems {\rmfamily \ref{T-1}} and {\rmfamily \ref{T-2}}.)
For any $n \geq 3$,
\begin{enumerate}
\item $\mathcal{E}_{F_n}(1) = \mathrm{Inn}\,F_n \cdot \mathcal{A}_{F_n}(2)$.
\item $\mathcal{A}_{F_n}(2k) \subset \mathcal{E}_{F_n}(k)$.
\end{enumerate}
\end{thrm}

From Proposition {\rmfamily \ref{I-P-2}}, we see that $\{ \mathcal{E}_G(k) \}$ is a central filtration of $\mathcal{E}_{G}(1)$.
Then a natural problem to consider is how different is $\{ \mathcal{E}_G(k) \}$ from the Andreadakis-Johnson filtration
$\{ \mathcal{A}_G(k) \}$. The partial answer to this question for $G=F_n$ is the theorem above.

\vspace{0.5em}

On the other hand, since $\{ \mathcal{E}_G(k) \}$ is central, each of the graded quotient
$\mathrm{gr}^k(\mathcal{E}_{F_n}):=\mathcal{E}_G(k)/\mathcal{E}_G(k+1)$ is an abelian group. At the end of the paper, we show
\begin{thrm}($=$ Theorem {\rmfamily \ref{P-4}}.)
For any $n \geq 3$,
\begin{enumerate}
\item Each of $\mathrm{gr}^k(\mathcal{E}_{F_n})$ is torsion-free.
\item $\mathrm{dim}_{\Q} (\mathrm{gr}^k(\mathcal{E}_{F_n}) \otimes_{\Z} \Q) < \infty$.
\end{enumerate}
\end{thrm}
To show this, we introduce Johnson homomorphism like homomorphisms $\eta_k$.
Observing Theorem {\rmfamily \ref{I-T-1}}, we see that $\mathrm{gr}^1(\mathcal{E}_{F_n})$ is finitely generated.
In general, however, it seems to be quite a difficult to determine the 
structure of $\mathrm{gr}^k(\mathcal{E}_{F_n})$ even the case where $k=1$.

\tableofcontents

\section{Notation and conventions}\label{S-Not}

\vspace{0.5em}

Throughout the paper, we use the following notation and conventions. Let $G$ be a group and $N$ a normal subgroup of $G$.
\begin{itemize}
\item The abelianization of $G$ is denoted by $G^{\mathrm{ab}}$.
\item The group $\mathrm{Aut}\,G$ of $G$ acts on $G$ from the right.
      For any $\sigma \in \mathrm{Aut}\,G$ and $x \in G$, the action of $\sigma$ on $x$ is denoted by $x^{\sigma}$.
\item For an element $g \in G$, we also denote the coset class of $g$ by $g \in G/N$ if there is no confusion.
\item For elements $x$ and $y$ of $G$, the commutator bracket $[x,y]$ of $x$ and $y$
      is defined to be $[x,y]:=xyx^{-1}y^{-1}$.
\end{itemize}
For pairs $(i_1, i_2, \ldots, i_k)$ and $(j_1, j_2, \ldots, j_k)$ of natural numbers $i_r, j_s \in \N$, we denote the lexicographic order
among them by $(i_1, i_2, \ldots, i_k) \leq (j_1, j_2, \ldots, j_k)$. Namely, this means $i_1 < j_1$, $i_1 = j_1$ and $i_2 < j_2$, or and so on.

\vspace{0.5em}

\section{The Andreadakis-Johnson filtration of $\mathrm{Aut}\,G$}\label{S-John}

\vspace{0.5em}

In this section, we review the Andreadakis-Johnson filtration of $\mathrm{Aut}\,G$ without proofs. The main purpose of the section is to fix the notations.
For basic materials concerning the Andreadakis-Johnson filtration and the Johnson homomorphisms, see \cite{S06} or \cite{S14}, for example.

\vspace{0.5em}

For a group $G$, we define the lower central series of $G$ by the rule
\[ \Gamma_G(1):= F_n, \hspace{1em} \Gamma_G(k) := [\Gamma_G(k-1),G], \hspace{1em} k \geq 2. \]
For any $y_1, \ldots, y_k \in G$, a left-normed commutator
\[ [[ \cdots [[y_1, y_2], y_3], \cdots ], y_k] \]
of weight $k$ is denoted by
\[ [y_1, y_2, \cdots, y_k ] \]
for simplicity. Then we have
\begin{lem}[See Section 5.3 in \cite{MKS}.]\label{l-1}
For any $k \geq 1$, the group $\Gamma_G(k)$ is generated by all left-normed commutators of weight $k$.
\end{lem}

\vspace{0.5em}

Let $\rho_G : \mathrm{Aut}\,G \rightarrow \mathrm{Aut}\,G^{\mathrm{ab}}$ be the natural homomorphism induced from
the abelianization of $G$. The kernel $\mathrm{IA}(G)$ of $\rho_G$ is called the IA-automorphism group of $G$.
Similarly, for any $k \geq 1$, the action of $\mathrm{Aut}\,G$ on each nilpotent quotient group $G/\Gamma_G(k+1)$ induces a homomorphism
\[ \rho_G^k : \mathrm{Aut}\,G \rightarrow \mathrm{Aut}(G/\Gamma_G(k+1)). \]
We denote the kernel of $\rho_G^k$ by $\mathcal{A}_G(k)$. Then the groups $\mathcal{A}_G(k)$ define a descending central filtration
\[ \mathrm{IA}(G) = \mathcal{A}_G(1) \supset \mathcal{A}_G(2) \supset \cdots \]
of $\mathrm{Aut}\,G$. We call it the Andreadakis-Johnson filtration of $\mathrm{Aut}\,G$. Then we have
\begin{thm}[Andreadakis \cite{And}]\label{T-And}
For any $k, l \geq 1$, $[\mathcal{A}_G(k), \mathcal{A}_G(l)] \subset \mathcal{A}_G(k+l)$.
\end{thm}

\vspace{0.5em}

\section{The ring of Fricke characters}\label{S-Fri}

\vspace{0.5em}

In this section, we review the ring of Fricke characters of a finitely generated group $G$.
In particular, we introduce a descending filtration consisting of $\mathrm{Aut}\,G$-invariant ideals of the ring.

\vspace{0.5em}

\subsection{The $\mathrm{Aut}\,G$-invariant ideal $J$}\label{Ss-Des}
\hspace*{\fill}\ 

\vspace{0.5em}

Let $G$ be a group generated by elements $x_1, \ldots, x_n$. We denote by
\[ R(G) := \mathrm{Hom}(G, \mathrm{SL}(2,\C)) \]
the set of group homomorphisms from $G$ to $\mathrm{SL}(2,\C)$. Let
\[ \mathcal{F}(R(G),\C) := \{ \chi : R(G) \rightarrow \C \} \]
be the set of complex-valued functions of $R(G)$. Then we can consider $\mathcal{F}(R(G),\C)$ as a $\C$-algebra in the following usual manner.
Namely, for any $\chi$, $\chi' \in \mathcal{F}(R(G),\C)$ and $\lambda \in \C$, we have
\[\begin{split}
  (\chi+ \chi')(\rho) & := \chi(\rho) + \chi'(\rho), \\
  (\chi \chi')(\rho) & := \chi(\rho) \chi'(\rho), \\
  (\lambda \chi)(\rho) & := \lambda \chi(\rho)
  \end{split}\]
for any $\rho \in R(G)$.
We consider $R(G)$ and $\mathcal{F}(R(G),\C)$ as right $\mathrm{Aut}(G)$-modules by
\[ \rho^{\sigma}(x) := \rho(x^{\sigma^{-1}}), \hspace{1em} \rho \in R(G) \,\,\, \text{and} \,\,\, x \in G \]
and
\[ \chi^{\sigma}(\rho) := \chi(\rho^{\sigma^{-1}}), \hspace{1em} \chi \in \mathcal{F}(R(G),\C) \,\,\, \text{and} \,\,\, \rho \in R(G) \]
respectively.

\vspace{0.5em}

For any $x \in G$, we define an element $\mathrm{tr}\,x$ of $\mathcal{F}(R(G),\C)$ to be
\[ (\mathrm{tr}\,x)(\rho) := \mathrm{tr}\, \rho(x) \]
for any $\rho \in R(G)$. Here \lq\lq $\mathrm{tr}$" in the right hand side means the trace of $2 \times 2$ matrix $\rho(x)$. The element $\mathrm{tr}\,x$ in
$\mathcal{F}(R(G),\C)$ is called the Fricke character of $x \in G$. 
The action of an element $\sigma \in \mathrm{Aut}(G)$ on $\mathrm{tr}\,x$ is given by $\mathrm{tr}\,x^{\sigma}$.
We have the following well-known formulae:
\begin{eqnarray}
 & & \mathrm{tr}\,x^{-1} = \mathrm{tr}\,x, \label{eq-1} \\
 & & \mathrm{tr}\,xy = \mathrm{tr}\,yx, \label{eq-2} \\
 & & \mathrm{tr}\,xy + \mathrm{tr}\,xy^{-1} =(\mathrm{tr}\,x)(\mathrm{tr}\,y), \label{eq-3} \\
 & & \mathrm{tr}\,xyz + \mathrm{tr}\,yxz =(\mathrm{tr}\,x)(\mathrm{tr}\,yz) + (\mathrm{tr}\,y)(\mathrm{tr}\,xz) 
      + (\mathrm{tr}\,z)(\mathrm{tr}\,xy) - (\mathrm{tr}\,x)(\mathrm{tr}\,y)(\mathrm{tr}\,z), \label{eq-4} \\
 & & \mathrm{tr}\,[x, y] = (\mathrm{tr}\,x)^2 + (\mathrm{tr}\,y)^2 + (\mathrm{tr}\,xy)^2 - (\mathrm{tr}\,x)(\mathrm{tr}\,y)(\mathrm{tr}\,xy) -2, \label{eq-4.5} \\
 & & 2\mathrm{tr}\,xyzw =(\mathrm{tr}\,x)(\mathrm{tr}\,yzw) + (\mathrm{tr}\,y)(\mathrm{tr}\,zwx)
       + (\mathrm{tr}\,z)(\mathrm{tr}\,wxy) + (\mathrm{tr}\,w)(\mathrm{tr}\,xyz) \label{eq-5} \\ 
 & & \hspace{5em} + (\mathrm{tr}\,xy)(\mathrm{tr}\,zw) - (\mathrm{tr}\,xz)(\mathrm{tr}\,yw) + (\mathrm{tr}\,xw)(\mathrm{tr}\,yz) \nonumber \\
 & & \hspace{5em} - (\mathrm{tr}\,x)(\mathrm{tr}\,y)(\mathrm{tr}\,zw) - (\mathrm{tr}\,y)(\mathrm{tr}\,z)(\mathrm{tr}\,xw)
                  - (\mathrm{tr}\,x)(\mathrm{tr}\,w)(\mathrm{tr}\,yz) \nonumber \\
 & & \hspace{5em} - (\mathrm{tr}\,z)(\mathrm{tr}\,w)(\mathrm{tr}\,xy) + (\mathrm{tr}\,x)(\mathrm{tr}\,y)(\mathrm{tr}\,z)(\mathrm{tr}\,w) \nonumber
\end{eqnarray}
for any $x, y, z, w \in G$. The equations (\ref{eq-4}) and (\ref{eq-5}) are due to Vogt \cite{Vog}.
(For details, see Section 3.4 in \cite{MaR} for example.) Note that
the whole point of (\ref{eq-4}) is that $\mathrm{tr}\,yxz$ can be written as a sum of $-\mathrm{tr}\,xyz$ and a polynomial in
$\mathrm{tr}\,v$ with $v$ a word in $x$, $y$, $z$ of length at most two. Similarly,
the whole point of (\ref{eq-5}) is that $\mathrm{tr}\,xyzw$ can be written as a polynomial in
$\mathrm{tr}\,v$ with $v$ a word in $x$, $y$, $z$, $w$ of length at most three.

\vspace{0.5em}

Let $\mathfrak{X}(G)$ be the $\Z$-submodule of $\mathcal{F}(R(G),\C)$ generated by all $\mathrm{tr}\,x$ for $x \in G$. Then, from
(\ref{eq-3}), it is seen that
$\mathfrak{X}(G)$ is closed under the multiplication of $\mathcal{F}(R(G),\C)$.
Consider an integral polynomial ring
\[ \Z[t] := \Z[t_{i_1 \cdots i_l} \,|\, 1 \leq l \leq n, \,\, 1 \leq i_1 < i_2 < \cdots < i_l \leq n] \]
of $2^n-1$ indeterminates, and the ring homomorphism $\pi : \Z[t] \rightarrow \mathcal{F}(R(G),\C)$ defined by
\[ \pi(1) := \frac{1}{2} (\mathrm{tr}\, 1_G), \hspace{1em} \pi(t_{i_1 \cdots i_l}) := \mathrm{tr}\,x_{i_1} \cdots x_{i_l}. \]
Fricke \cite{Fri} stated that for any element $x$ which is a word $x \in G$ in the generators $x_1, \ldots x_n$, the character $\mathrm{tr}\,x$
is a polynomial among $\mathrm{tr}\,x_{i_1} x_{i_2} \cdots x_{i_l}$ for $1 \leq l \leq n$ and $1 \leq i_1 < i_2 < \cdots < i_l \leq n$.
This was proved by Horowitz \cite{Ho1}. More precisely,
\begin{thm}[Horowitz, \cite{Ho1}]\label{T-Hor1}
For any $G$, $\mathfrak{X}(G)$ is the image of an ideal
\[ I_0 := (2, t_{i_1 \cdots i_l} \,|\, 1 \leq l \leq n, \,\, 1 \leq i_1 < i_2 < \cdots < i_l \leq n) \subset \Z[t] \]
by $\pi$. 
\end{thm}
Set $I := \mathrm{Ker}(\pi)$. Namely,
\[ I = \{ f \in \Z[t] \,|\, f(\mathrm{tr}\,\rho(x_{i_1} \cdots x_{i_l}))=0 \,\,\, \text{for any} \,\,\, \rho \in R(G) \}. \]
Consider the case where $G=F_n$. Horowitz \cite{Ho1} also showed that $I=(0)$ for $n=1$ and $2$, and that $I$ is a principal ideal generated by a quadratic element
\[ t_{123}^2 - P_{123}(t) t_{123} + Q_{123}(t) \]
where
\[\begin{split}
  P_{abc}(t) & := t_{ab} t_c + t_{ac} t_b + t_{bc} t_a, \\
  Q_{abc}(t) & := t_a^2 + t_b^2 + t_c^2 + t_{ab}^2 + t_{ac}^2 + t_{bc}^2
                   - t_a t_b t_{ab} - t_a t_c t_{ac} - t_b t_c t_{bc} + t_{ab} t_{bc} t_{ac} -4
 \end{split}\]
for $n=3$.
For $n \geq 4$, Whittemore \cite{Whi} showed that $I$ is not principal. In general, however,
very little is known for the ideal $I$ for general $n \geq 4$.

\vspace{0.5em}

In this paper, we call the quotient ring $\Z[t]/I$ the ring of Fricke characters of $G$ over $\Z$, and
considered as a subring of $\mathcal{F}(R(G),\C)$ through the homomorphism $\pi$.
Then, we can define the $\mathrm{Aut}\,(G)$-module structure of $\Z[t]/I$ such that the induced homomorphism
$\Z[t]/I \rightarrow \mathcal{F}(R(G),\C)$ from $\pi$ is $\mathrm{Aut}\,(G)$-equivariant injective.

\vspace{0.5em}

For an elements $y \in G$, an automorphism $\iota_y$ of $G$ defined by $x^{\iota_y} := y^{-1} xy$ for any $x \in G$ is called an inner automorphism of $G$
associated to $y$.
Let $\mathrm{Inn}\,(G)$ be a normal subgroup of $\mathrm{Aut}\,(G)$ consisting of inner automorphisms of $G$.
In general, $\mathrm{Inn}\,(G)$ is contained in the kernel of the
homomorphism $\zeta : \mathrm{Aut}\,G \rightarrow \mathrm{Aut}(\Z[t]/I)$
induced from the action of $\mathrm{Aut}\,G$ on the ring of Fricke characters.
For the case where $G=F_n$, Horowitz \cite{Ho2} showed
\begin{thm}[Horowitz, \cite{Ho2}]\label{T-Hor2}
For $n \geq 3$, $\mathrm{Ker}(\zeta) =\mathrm{Inn}\,F_n$.
\end{thm}
Namely, the action of $\mathrm{Aut}\,F_n$ on the ring of Fricke characters induces a faithful representation of the outer automorphism group
$\mathrm{Out}\,F_n := \mathrm{Aut}\,F_n/\mathrm{Inn}\,F_n$ of a free group $F_n$.
However, since $\Z[t]/I$ in not finitely generated as a $\Z$-module, the representations
$\mathrm{Aut}\,F_n \rightarrow \mathrm{Aut}(\Z[t]/I)$ and $\mathrm{Out}\,F_n \rightarrow \mathrm{Aut}(\Z[t]/I)$ are not so easy to handle in general.
In addition to this, the number of indeterminates of $\Z[t]$ also adds momentum to the complexity if we write down the behavior of the actions of
$\mathrm{Aut}\,F_n$ and $\mathrm{Out}\,F_n$ on $\Z[t]/I$.

\vspace{0.5em}

In order to avoid these difficulties, first, we consider the rationalization of the situation above.
Let $\mathfrak{X}_{\Q}(G)$ be the $\Q$-subspace of $\mathcal{F}(R(G),\C)$ generated by all $\mathrm{tr}\,x$ for $x \in G$.
The set $\mathfrak{X}_{\Q}(G)$ naturally has a ring structure.
Let
\[ \Q[t] := \Q[t_{i_1 \cdots i_l} \,|\, 1 \leq l \leq 3, \,\, 1 \leq i_1 < i_2 < \cdots < i_l \leq n] \]
be a rational polynomial ring of $n + \binom{n}{2} + \binom{n}{3}$ indeterminates.
Consider the ring homomorphism $\pi_{\Q} : \Q[t] \rightarrow \mathcal{F}(R(G),\C)$ defined by
\[ \pi_{\Q}(1) := \frac{1}{2} (\mathrm{tr}\, 1_G), \hspace{1em} \pi_{\Q}(t_{i_1 \cdots i_l}) := \mathrm{tr}\,x_{i_1} \cdots x_{i_l}. \]
Then observing the formula (\ref{eq-5}) and Horowitz's result as mentioned above, we see $\mathrm{Im}(\pi_{\Q}) = \mathfrak{X}_{\Q}(G)$.
Remark that $\mathrm{Im}(\pi) \neq \mathfrak{X}(G)$.
Set
\[ I_{\Q} := \mathrm{Ker}(\pi_{\Q})
   = \{ f \in \Q[t] \,|\, f(\mathrm{tr}\,\rho(x_{i_1} \cdots x_{i_l}))=0 \,\,\, \text{for any} \,\,\, \rho \in R(G) \}. \]
We call $\mathfrak{X}_{\Q}(G)$ and $\Q[t]/I_{\Q}$ the ring of Fricke characters of $G$ over $\Q$.
Similar to $\Z[t]/I$, we see that $\Q[t]/I_{\Q}$ can be considered as an $\mathrm{Aut}\,(G)$-module, and that
$\mathrm{Inn}\,(G)$ is contained in the kernel of the homomorphism $\zeta_{\Q} : \mathrm{Aut}\,G \rightarrow \mathrm{Aut}(\Q[t]/I)$
induced from the action of $\mathrm{Aut}\,G$ on $\Q[t]/I$.

\vspace{0.5em}

If $G=F_n$, since $\mathrm{Ker}(\zeta_{\Q})$ acts on $\mathfrak{X}(F_n) \subset \mathfrak{X}_{\Q}(F_n)$ trivially,
we see that $\mathrm{Ker}(\zeta_{\Q})=\mathrm{Inn}\,F_n$ by Theorem {\rmfamily \ref{T-Hor2}}.
Hence, $\zeta_{\Q}$ also induces a faithful representation of $\mathrm{Out}\,F_n$.
In order to construct finite dimensional representations of $\mathrm{Aut}\,G$ and $\mathrm{Out}\,G$, we prepare a descending filtration of 
$\mathrm{Aut}\,G$-invariant ideals of $\Q[t]/I_{\Q}$, and take its graded quotients.
Set $t_{i_1 \cdots i_l}' := t_{i_1 \cdots i_l} -2 \in \Q[t]$. We also denote by $t_{i_1 \cdots i_l}'$ its coset class in $\Q[t]/I_{\Q}$ for simplicity.
Consider the ideal
\[ J := (t_{i_1 \cdots i_l}' \,|\, 1 \leq l \leq 3, \,\, 1 \leq i_1 < i_2 < \cdots < i_l \leq n) \subset \Q[t]/I_{\Q} \]
generated by all $t_{i_1 \cdots i_l}'$'s.
\begin{lem}[For $n=3$, see also Magnus \cite{Mg2}.]
The ideal $J$ is $\mathrm{Aut}\,G$-invariant.
\end{lem}
\textit{Proof.}
For any $t_{j_1 \cdots j_m}'$ and $\sigma \in \mathrm{Aut}\,G$, there exists some polynomial $F(t_{i_1 \cdots i_l}) \in \Q[t]$ such that
\[ (t_{j_1 \cdots j_m}')^{\sigma} \equiv F(t_{i_1 \cdots i_l}) \in \Q[t]/I_{\Q} \]
by Theorem {\rmfamily \ref{T-Hor1}} and (\ref{eq-5}). Then using the division algorithm, we verify that $F$ can be written as
\[ F = t_1' G_1 + R \in \Q[t] \]
where $G_1$, $R \in \Q[t]$ such that $R$ is a polynomial in the determinates $t_{i_1 \cdots i_l}'$ except for $t_1'$.
By repeating this argument, we obtain
\[ F = \sum_{i_1< \cdots <i_l} t_{i_1 \cdots i_l}' G_{i_1 \cdots i_l} + C \]
where $G_{i_1 \cdots i_l} \in \Q[t]$, $C \in \Q$, and the sum runs over all $i_1< \cdots <i_l$ such that $1 \leq l \leq 3$ and $1 \leq i_1 < i_2 < \cdots < i_l \leq n$.

\vspace{0.5em}

By considering the image of this equation by $\zeta_{\Q}$, and by substituting the trivial representation
${\bf 1} : R(G) \rightarrow \mathrm{SL}(2, \C)$, we see that $C=0$. This means $F \in J$. Therefore, $J$ is $\mathrm{Aut}\,G$-invariant. $\square$

\vspace{0.5em}

Now, observe the descending filtration
\[ J \supset J^2 \supset J^3 \supset \cdots \]
of $\Q[t]/I_{\Q}$, and set $\mathrm{gr}^k(J) := J^k/J^{k+1}$ for any $k \geq 1$.
Then each of $\mathrm{gr}^k(J)$ is an $\mathrm{Aut}\,G$-invariant $\Q$-vector space of finite dimension.
Hence, we obtain finite dimensional representations
\[ \zeta_{k,\Q} : \mathrm{Aut}\,G \rightarrow \mathrm{Aut}\,(\mathrm{gr}^k(J)) \]
over $\Q$.

\vspace{0.5em}

Consider the case where $G=F_n$. First,
since $\mathrm{Ker}(\zeta_{\Q})$ acts on $\mathfrak{X}(F_n) \subset \mathfrak{X}_{\Q}(F_n)$ trivially,
we see that $\mathrm{Ker}(\zeta_{\Q})=\mathrm{Inn}\,F_n$ by Theorem {\rmfamily \ref{T-Hor2}}.
Hence, $\zeta_{\Q}$ also induces a faithful representation of $\mathrm{Out}\,F_n$.
On the other hand, in order to investigate the behavior of the action of $\mathrm{Aut}\,F_n$ on $\mathrm{gr}^k(J)$,
we have to know the $\Q$-vector space structures of $\mathrm{gr}^k(J)$.
By combinatorial complexities, however, it seems quite difficult to obtain a basis of $\mathrm{gr}^k(J)$ for a general $k$.
For $k=1$ and $2$, we give bases of $\mathrm{gr}^k(J)$ in Subsection {\rmfamily \ref{Ss-Str}}.

\vspace{0.5em}

\subsection{Basic formulae among $\mathrm{tr}'\,x$}\label{Ss-For}
\hspace*{\fill}\ 

\vspace{0.5em}

For any $x \in G$, set
\[ \mathrm{tr}'\,x := (\mathrm{tr}\,x) -2 \in \mathcal{F}(R(G),\C). \]
In this subsection, we summerize basic and useful formulae among $\mathrm{tr}'\,x$. To begin with, in order to rewrite
the equations (\ref{eq-1}), $\ldots$, (\ref{eq-5}) as those among $\mathrm{tr}'\,x, \ldots, \mathrm{tr}'\,w$, we prepare
\begin{lem}
For any $k \geq 1$, and $z_1, \ldots, z_k \in G$, we have
\[\begin{split}
   (\mathrm{tr}\,z_1) \cdots (\mathrm{tr}\,z_k)
     = \sum_{i=0}^k 2^i \sum_{1 \leq j_1 < \cdots < j_{k-i} \leq k} (\mathrm{tr}'\,z_{j_1}) \cdots (\mathrm{tr}'\,z_{j_{k-i}})
  \end{split}\]
\end{lem}
Since this formula can be shown easily with the induction on $k$, we omit the details. $\square$ 

\vspace{0.5em}

Then using the lemma above and (\ref{eq-1}), $\ldots$, (\ref{eq-5}), we obtain 
\begin{eqnarray}
 & & \mathrm{tr}'\,x^{-1} = \mathrm{tr}'\,x, \label{eq-6} \\
 & & \mathrm{tr}'\,xy = \mathrm{tr}'\,yx, \label{eq-7} \\
 & & \mathrm{tr}'\,xy + \mathrm{tr}'\,xy^{-1} =2\mathrm{tr}'\,x + 2\mathrm{tr}'\,y + (\mathrm{tr}'\,x)(\mathrm{tr}'\,y), \label{eq-8} \\
 & & \mathrm{tr}'\,xyz + \mathrm{tr}'\,yxz = -2 \{ \mathrm{tr}'\,x + \mathrm{tr}'\,y + \mathrm{tr}'\,z \}
                                                  + 2 \{ \mathrm{tr}'\,xy + \mathrm{tr}'\,yz + \mathrm{tr}'\,xz \} \label{eq-9} \\
 & & \hspace{8em} + (\mathrm{tr}'\,x)(\mathrm{tr}'\,yz) + (\mathrm{tr}'\,y)(\mathrm{tr}'\,xz) 
      + (\mathrm{tr}'\,z)(\mathrm{tr}'\,xy), \nonumber \\
 & & \hspace{8em} -2 \{ (\mathrm{tr}'\,x)(\mathrm{tr}'\,y) + (\mathrm{tr}'\,y)(\mathrm{tr}'\,z) + (\mathrm{tr}'\,z)(\mathrm{tr}'\,x) \} \nonumber \\
 & & \hspace{8em} - (\mathrm{tr}'\,x)(\mathrm{tr}'\,y)(\mathrm{tr}'\,z), \nonumber \\
 & & \mathrm{tr}'\,[x, y] = (\mathrm{tr}'\,x)^2 + (\mathrm{tr}'\,y)^2 + (\mathrm{tr}'\,xy)^2 \label{eq-10} \\
 & & \hspace{5em}  -2 \{ (\mathrm{tr}'\,x)(\mathrm{tr}'\,y) + (\mathrm{tr}'\,x)(\mathrm{tr}'\,xy) + (\mathrm{tr}'\,y)(\mathrm{tr}'\,xy) \} 
        - (\mathrm{tr}'\,x)(\mathrm{tr}'\,y)(\mathrm{tr}'\,xy)  \nonumber
\end{eqnarray}
and
{\small
\begin{equation}\begin{split}
2\mathrm{tr}'\,xyzw & = 2(\mathrm{tr}'\,x + \mathrm{tr}'\,y + \mathrm{tr}'\,z + \mathrm{tr}'\,w) \\
   & \hspace{2em} - 2(\mathrm{tr}'\,xy + \mathrm{tr}'\,xz + \mathrm{tr}'\,xw + \mathrm{tr}'\,yz + \mathrm{tr}'\,yw + \mathrm{tr}'\,zw) \\
   & \hspace{2em} +2(\mathrm{tr}'\,xyz + \mathrm{tr}'\,xyw + \mathrm{tr}'\,xzw + \mathrm{tr}'\,yzw) \\
   & \hspace{2em} +2 \{ (\mathrm{tr}'\,x)(\mathrm{tr}'\,y) + (\mathrm{tr}'\,x)(\mathrm{tr}'\,w) + (\mathrm{tr}'\,y)(\mathrm{tr}'\,z)
                     + (\mathrm{tr}'\,z)(\mathrm{tr}'\,w) \\
   & \hspace{5em} + 2 (\mathrm{tr}'\,x)(\mathrm{tr}'\,z) + 2 (\mathrm{tr}'\,y)(\mathrm{tr}'\,w) \} \\
   & \hspace{2em} - 2 \{ (\mathrm{tr}'\,x)(\mathrm{tr}'\,yz) + (\mathrm{tr}'\,x)(\mathrm{tr}'\,zw) + (\mathrm{tr}'\,y)(\mathrm{tr}'\,xw)
                      + (\mathrm{tr}'\,y)(\mathrm{tr}'\,zw) \\
   & \hspace{5em} + (\mathrm{tr}'\,z)(\mathrm{tr}'\,xy) + (\mathrm{tr}'\,z)(\mathrm{tr}'\,xw) + (\mathrm{tr}'\,w)(\mathrm{tr}'\,xy)
                      + (\mathrm{tr}'\,w)(\mathrm{tr}'\,yz) \} \\
   & \hspace{2em} + \{ (\mathrm{tr}'\,x)(\mathrm{tr}'\,yzw) + (\mathrm{tr}'\,y)(\mathrm{tr}'\,xzw) + (\mathrm{tr}'\,z)(\mathrm{tr}'\,xyw)
                      + (\mathrm{tr}'\,w)(\mathrm{tr}'\,xyz) \} \\
   & \hspace{2em} + \{ (\mathrm{tr}'\,xy)(\mathrm{tr}'\,zw) - (\mathrm{tr}'\,xz)(\mathrm{tr}'\,yw) + (\mathrm{tr}'\,xw)(\mathrm{tr}'\,yz) \} \\
   & \hspace{2em} - \{ (\mathrm{tr}'\,x)(\mathrm{tr}'\,y)(\mathrm{tr}'\,zw) + (\mathrm{tr}'\,y)(\mathrm{tr}'\,z)(\mathrm{tr}'\,xw) 
                       + (\mathrm{tr}'\,x)(\mathrm{tr}'\,w)(\mathrm{tr}'\,yz) \\
   & \hspace{5em} + (\mathrm{tr}'\,z)(\mathrm{tr}'\,w)(\mathrm{tr}'\,xy) \} \\
   & \hspace{2em} + (\mathrm{tr}'\,x)(\mathrm{tr}'\,y)(\mathrm{tr}'\,z)(\mathrm{tr}'\,w) \\
   & \hspace{2em} + 2 \{ (\mathrm{tr}'\,x)(\mathrm{tr}'\,y)(\mathrm{tr}'\,z) + (\mathrm{tr}'\,x)(\mathrm{tr}'\,y)(\mathrm{tr}'\,w)
                  + (\mathrm{tr}'\,x)(\mathrm{tr}'\,z)(\mathrm{tr}'\,w) \\
   & \hspace{5em} + (\mathrm{tr}'\,y)(\mathrm{tr}'\,z)(\mathrm{tr}'\,w) \}.
\end{split} \label{eq-20} \end{equation}
}

\vspace{0.5em}

Furthermore, we can rewrite (\ref{eq-20}) as a sum of the degree one part and elements types of
\[ (\mathrm{tr}'\,\alpha)(\mathrm{tr}'\,\beta w - \mathrm{tr}'\,\beta), \,\,\, (\mathrm{tr}'\,\alpha)(\mathrm{tr}'\,w), \,\,\,
   (\mathrm{tr}'\,\alpha)(\mathrm{tr}'\,\beta)(\mathrm{tr}'\,w) \]
for some $\alpha, \beta \in G$ and
\[ (\mathrm{tr}'\,x)(\mathrm{tr}'\,y)(\mathrm{tr}'\,z)(\mathrm{tr}'\,w). \]
That is,
{\small
\begin{equation}\begin{split}
2\mathrm{tr}'& \,xyzw = 2(\mathrm{tr}'\,x + \mathrm{tr}'\,y + \mathrm{tr}'\,z + \mathrm{tr}'\,w) \\
   & \hspace{2em} - 2(\mathrm{tr}'\,xy + \mathrm{tr}'\,xz + \mathrm{tr}'\,xw + \mathrm{tr}'\,yz + \mathrm{tr}'\,yw + \mathrm{tr}'\,zw) \\
   & \hspace{2em} +2(\mathrm{tr}'\,xyz + \mathrm{tr}'\,xyw + \mathrm{tr}'\,xzw + \mathrm{tr}'\,yzw) \\
   & \hspace{2em} +2 \{ (\mathrm{tr}'\,x - \mathrm{tr}'\,xw)(\mathrm{tr}'\,y) + (\mathrm{tr}'\,y)(\mathrm{tr}'\,z - \mathrm{tr}'\,zw)
                    + (\mathrm{tr}'\,x - \mathrm{tr}'\,xw)(\mathrm{tr}'\,z)  \\
   & \hspace{3em} + (\mathrm{tr}'\,x)(\mathrm{tr}'\,z - \mathrm{tr}'\,zw) \} \\
   & \hspace{2em} - (\mathrm{tr}'\,x)(\mathrm{tr}'\,yz - \mathrm{tr}'\,yzw) - (\mathrm{tr}'\,x- \mathrm{tr}'\,xw)(\mathrm{tr}'\,yz)
                  - (\mathrm{tr}'\,z - \mathrm{tr}'\,zw)(\mathrm{tr}'\,xy)   \\
   & \hspace{3em} - (\mathrm{tr}'\,z)(\mathrm{tr}'\,xy- \mathrm{tr}'\,xyw) 
                  + (\mathrm{tr}'\,y)(\mathrm{tr}'\,xzw - \mathrm{tr}'\,xz) + (\mathrm{tr}'\,xz)(\mathrm{tr}'\,y - \mathrm{tr}'\,yw) \\
   & \hspace{2em} + (\mathrm{tr}'\,x)(\mathrm{tr}'\,y)(\mathrm{tr}'\,z - \mathrm{tr}'\,zw)
                  + (\mathrm{tr}'\,y)(\mathrm{tr}'\,z)(\mathrm{tr}'\,x - \mathrm{tr}'\,xw) \\ \\
   & \hspace{2em} + 2(\mathrm{tr}'\,x)(\mathrm{tr}'\,w) + 2(\mathrm{tr}'\,z)(\mathrm{tr}'\,w)  + 4 (\mathrm{tr}'\,y)(\mathrm{tr}'\,w)
     - 2(\mathrm{tr}'\,w)(\mathrm{tr}'\,xy)  \\
   & \hspace{2em} - 2 (\mathrm{tr}'\,w)(\mathrm{tr}'\,yz) + (\mathrm{tr}'\,w)(\mathrm{tr}'\,xyz) \\
   & \hspace{2em} + (\mathrm{tr}'\,x)(\mathrm{tr}'\,y)(\mathrm{tr}'\,z)(\mathrm{tr}'\,w) + (\mathrm{tr}'\,x)(\mathrm{tr}'\,w)(\mathrm{tr}'\,yz)
                  + (\mathrm{tr}'\,z)(\mathrm{tr}'\,w)(\mathrm{tr}'\,xy)  \\
   & \hspace{2em} + 2 \{ (\mathrm{tr}'\,x)(\mathrm{tr}'\,y)(\mathrm{tr}'\,w) + (\mathrm{tr}'\,x)(\mathrm{tr}'\,z)(\mathrm{tr}'\,w)
                  + (\mathrm{tr}'\,y)(\mathrm{tr}'\,z)(\mathrm{tr}'\,w) \}.
\end{split} \label{eq-30} \end{equation}
}
\vspace{0.5em}

Next, we consider elements type of $\mathrm{tr}' z \gamma$ for $z \in G$ and $\gamma \in [G,G]$. First, we study the case where $\gamma$ is a
commutator of weight $2$.
\begin{lem}\label{L-9}
For any $z, a, b \in G$,
\[\begin{split}
  \mathrm{tr}' z[a,b] = & \,\, \mathrm{tr}' z + 2(\mathrm{tr}' z + \mathrm{tr}' a + \mathrm{tr}' b) \\
    & -2(\mathrm{tr}' za + \mathrm{tr}' zb + \mathrm{tr}' ab) + 2 \mathrm{tr}' zab \\
    & +(\mathrm{tr}' za)(\mathrm{tr}' a) - (\mathrm{tr}' zb)(\mathrm{tr}' b) + 4 (\mathrm{tr}' z)(\mathrm{tr}' b) +2 (\mathrm{tr}' b)^2 \\
    & -2 (\mathrm{tr}' za)(\mathrm{tr}' b) -2 (\mathrm{tr}' ab)(\mathrm{tr}' b) -2 (\mathrm{tr}' za)(\mathrm{tr}' ab) \\
    & + (\mathrm{tr}' ab)(\mathrm{tr}' zab) + (\mathrm{tr}' z)(\mathrm{tr}' b)^2 - (\mathrm{tr}' za)(\mathrm{tr}' ab)(\mathrm{tr}' b) \\
  = & \,\, \mathrm{tr}' z -2(\mathrm{tr}' za - \mathrm{tr}' z) +2(\mathrm{tr}' bza - \mathrm{tr}' bz) -2(\mathrm{tr}' ba - \mathrm{tr}' b) \\
   & + (\mathrm{tr}' bza - \mathrm{tr}' bz)(\mathrm{tr}' b) + (\mathrm{tr}' ba - \mathrm{tr}' b)(\mathrm{tr}' bza)
     - 2(\mathrm{tr}' za - \mathrm{tr}' z)(\mathrm{tr}' b) \\
   & - 2(\mathrm{tr}' ba - \mathrm{tr}' b)(\mathrm{tr}' b) - 2(\mathrm{tr}' za - \mathrm{tr}' z)(\mathrm{tr}' b) -2 (\mathrm{tr}' ba - \mathrm{tr}' b)(\mathrm{tr}' za) \\
   & - (\mathrm{tr}' za - \mathrm{tr}' z)(\mathrm{tr}' b)^2 - (\mathrm{tr}' ba - \mathrm{tr}' b)(\mathrm{tr}' b)(\mathrm{tr}' za) \\
   & + 2 \mathrm{tr}' a + (\mathrm{tr}' za)(\mathrm{tr}' a).
 \end{split}\]
\end{lem}

\textit{Proof.}
We show the former equality. The latter one immediately follows from the former one. 
Now, we have
\[\begin{split}
   \mathrm{tr}' z[a,b] & = \mathrm{tr}' zab(ba)^{-1} \\
    & \stackrel{(\ref{eq-8})}{=} - \mathrm{tr}' zab^2 a + (\mathrm{tr}' zab)(\mathrm{tr}' ab) + 2 \mathrm{tr}' zab +2 \mathrm{tr}' ab \\
    & = - \mathrm{tr}' (azab)b + (\mathrm{tr}' zab)(\mathrm{tr}' ab) + 2 \mathrm{tr}' zab +2 \mathrm{tr}' ab \\
    & \stackrel{(\ref{eq-8})}{=} \mathrm{tr}' aza - (\mathrm{tr}' azab)(\mathrm{tr}' b) - 2 \mathrm{tr}' azab -2 \mathrm{tr}' b \\
    & \hspace{2em} + (\mathrm{tr}' zab)(\mathrm{tr}' ab) + 2 \mathrm{tr}' zab +2 \mathrm{tr}' ab \\
    & = \mathrm{tr}' (za)a - (\mathrm{tr}' azab)(\mathrm{tr}' b) - 2 \mathrm{tr}' azab -2 \mathrm{tr}' b \\
    & \hspace{2em} + (\mathrm{tr}' zab)(\mathrm{tr}' ab) + 2 \mathrm{tr}' zab +2 \mathrm{tr}' ab \\
    & \stackrel{(\ref{eq-8})}{=} - \mathrm{tr}' z + (\mathrm{tr}' za)(\mathrm{tr}' a) +2 \mathrm{tr}' za + 2 \mathrm{tr}' a \\
   & \hspace{2em} - \{ -\mathrm{tr}' azb^{-1}a^{-1} + (\mathrm{tr}' za)(\mathrm{tr}' ab) + 2 \mathrm{tr}' za +2 \mathrm{tr}' ab  \} (\mathrm{tr}' b) \\   &  \hspace{2em} - 2 \{ -\mathrm{tr}' azb^{-1}a^{-1} + (\mathrm{tr}' za)(\mathrm{tr}' ab) + 2 \mathrm{tr}' za + 2 \mathrm{tr}' ab  \} \\
   &  \hspace{2em} -2 \mathrm{tr}' b  + (\mathrm{tr}' zab)(\mathrm{tr}' ab) + 2 \mathrm{tr}' zab +2 \mathrm{tr}' ab.
  \end{split}\]
Using this and
\[ \mathrm{tr}' zb^{-1} \stackrel{(\ref{eq-8})}{=} - \mathrm{tr}' zb + (\mathrm{tr}' z)(\mathrm{tr}' b) +2 \mathrm{tr}' z + 2 \mathrm{tr}' b, \]
we obtain the required result. This completes the proof of Lemma {\rmfamily \ref{L-9}}. $\square$

\vspace{0.5em}

Using Lemma {\rmfamily \ref{L-9}} and (\ref{eq-9}), we see
\begin{cor}\label{C-9}
For any $z, a, b \in G$,
$\mathrm{tr}'\,z[a, b] \equiv \mathrm{tr}'\,z + \mathrm{tr}'\,zab - \mathrm{tr}'\,zba \pmod{J^2}$. 
\end{cor}

\vspace{0.5em}

Now, we consider the case where $\gamma$ is a commutator of weight $3$.
\begin{lem}\label{L-10}
For any $z, a, b \in G$, an element
$\mathrm{tr}'\,z[a, b, c] - \mathrm{tr}'\,z$ is a sum of elements types of
\[ (\mathrm{tr}'\,x)(\mathrm{tr}'\,y[a, b] - \mathrm{tr}'\,y), \,\,\, (\mathrm{tr}'\,x)(\mathrm{tr}'\,ya - \mathrm{tr}'\,y) \]
and
\[ \mathrm{tr}'\,[a, b], \,\,\, (\mathrm{tr}'\,x)(\mathrm{tr}'\,[a, b]), \,\,\, (\mathrm{tr}'\,x)(\mathrm{tr}'\,a) \]
for some $x, y \in G$.
\end{lem}
\textit{Proof.}
By substituting $a$ and $b$ in the equation in Lemma {\rmfamily \ref{L-9}} to $[a,b]$ and $c$ respectively, we obtain
\[\begin{split}
   \mathrm{tr}' z[a,b,c] = & \,\, \mathrm{tr}' z -2(\mathrm{tr}' z[a,b] - \mathrm{tr}' z) +2(\mathrm{tr}' cz[a,b] - \mathrm{tr}' cz)
         -2(\mathrm{tr}' c[a,b] - \mathrm{tr}' c) \\
   & + 2 \mathrm{tr}' [a,b] + (\mathrm{tr}' z[a,b])(\mathrm{tr}' [a,b]) + \sum (\mathrm{tr}' x)(\mathrm{tr}' y[a,b] - \mathrm{tr}\, y).
  \end{split}\]
Again, by Lemma {\rmfamily \ref{L-9}}, we have
\[\begin{split}
   \mathrm{tr}' z[a,b,c] = & \,\, \mathrm{tr}' z -2 \{ -2 (\mathrm{tr}' za - \mathrm{tr}' z) +2(\mathrm{tr}' bza - \mathrm{tr}' bz) -2(\mathrm{tr}' ba - \mathrm{tr}' b) \} \\
    & +2 \{ -2 (\mathrm{tr}' cza - \mathrm{tr}' cz) +2(\mathrm{tr}' bcza - \mathrm{tr}' bcz) -2(\mathrm{tr}' ba - \mathrm{tr}' b) \} \\
    & -2 \{ -2 (\mathrm{tr}' ca - \mathrm{tr}' c) +2(\mathrm{tr}' bca - \mathrm{tr}' bc) -2(\mathrm{tr}' ba - \mathrm{tr}' b) \} \\
    & -4 \mathrm{tr}' a + \sum (\mathrm{tr}' x)(\mathrm{tr}' a) + \sum (\mathrm{tr}' x)(\mathrm{tr}' ya - \mathrm{tr}\, y) \\
    & + 2 \mathrm{tr}' [a,b] + (\mathrm{tr}' z[a,b])(\mathrm{tr}' [a,b]) + \sum (\mathrm{tr}' x)(\mathrm{tr}' y[a,b] - \mathrm{tr}\, y).
  \end{split}\]
Then by using (\ref{eq-30}) for $\mathrm{tr}' bcza$, we obatain the required result.
This completes the proof of Lemma {\rmfamily \ref{L-10}}. $\square$

\vspace{0.5em}

In particular, we have
\begin{cor}\label{C-10}
$\mathrm{tr}'\,z[a, b, c]^{\pm1} \equiv \mathrm{tr}'\,z \pmod{J^2}$.
\end{cor}

\textit{Proof.}
For $\mathrm{tr}'\,z[a, b, c]^{-1} \equiv \mathrm{tr}'\,z \pmod{J^2}$, it suffices to consider
\[ \mathrm{tr}'\,z[a, b, c]^{-1}=\mathrm{tr}'\,z[cac^{-1}, cbc^{-1}, c^{-1}]. \]
This completes the proof of Corollary {\rmfamily \ref{C-10}}. $\square$

\vspace{0.5em}

The next proposition is the goal of this subsection. For any $k \geq 2$ and $y_1, \ldots, y_k \in G$, set
\[ a_k :=[ y_1, y_2, \ldots, y_k] \in \Gamma_G(k). \]
\begin{pro}\label{P-7}
With the notation above, we have
\begin{enumerate}
\item For any $k \geq 2$, $\mathrm{tr}'\,a_k \in J^{k-1}$.
\item For any $l \geq 2$ and any $b \in G$,
\[ \mathrm{tr}'\,ba_{2l-1} - \mathrm{tr}'\,b, \,\,\, \mathrm{tr}'\,ba_{2l} - \mathrm{tr}'\,b \in J^l. \]
\end{enumerate}
\end{pro}
\textit{Proof.}
We prove this proposition by the induction on $k$ and $l$. To begin with, for $k=l=2$, Part (1) and Part (2) follows from
(\ref{eq-10}) and Corollary {\rmfamily \ref{C-10}} respectively.
(Remark that if we set $y_1' := [y_1, y_2]$, then we can see $a_4 = a_3(y_1', y_3, y_4)$, and hence $\mathrm{tr}'\,ba_4 - \mathrm{tr}'\,b \in J^2$.)
Furthermore, Part (1) also holds for $k=3$ by (\ref{eq-10}).

\vspace{0.5em}

Assume $l \geq 2$ and $k=2l-2 \geq 2$, and assume that Part (1) is true for $k$ and $k+1$, and that the part (2) is true for $l$.
We show that
for any $k'$ such that $k+2 \leq k' \leq k+3$, Part (1) holds. By the inductive hypothesis, we have
\[\begin{split}
   \mathrm{tr}'\, a_{k'} & = \mathrm{tr}'\,[a_{k'-1}, y_{k'}] \\
    & \stackrel{(\ref{eq-10})}{=} (\mathrm{tr}'\,a_{k'-1})^2 + (\mathrm{tr}'\,y_{k'})^2 + (\mathrm{tr}'\,a_{k'-1}y_{k'})^2 \\
    & \hspace{2em} -2 \{ (\mathrm{tr}'\,a_{k'-1})(\mathrm{tr}'\,y_{k'}) + (\mathrm{tr}'\,a_{k'-1})(\mathrm{tr}'\,a_{k'-1}y_{k'}) 
     + (\mathrm{tr}'\,y_{k'})(\mathrm{tr}'\,a_{k'-1}y_{k'}) \} \\ 
    & \hspace{2em} - (\mathrm{tr}'\,a_{k'-1})(\mathrm{tr}'\,y_{k'})(\mathrm{tr}'\,a_{k'-1}y_{k'}) \\
    & \equiv (\mathrm{tr}'\,a_{k'-1}y_{k'} - \mathrm{tr}'\,y_{k'})^2 \pmod{J^{k'-1}} \\ 
    & \equiv 0 \pmod{J^{k'-1}}.
  \end{split}\]
This shows that for $k'=k+1, k+2$, Part (1) holds.

\vspace{0.5em}

Next, for $l'=l+1$, by Lemma {\rmfamily \ref{L-10}}, we can write
\[\begin{split}
   \mathrm{tr}' & \,ba_{2l'-1} - \mathrm{tr}'\,b \\
    & = \sum (\mathrm{tr}'\,x)(\mathrm{tr}'\,ya_{2l'-3} - \mathrm{tr}'\,y) + \sum (\mathrm{tr}'\,x)(\mathrm{tr}'\,ya_{2l'-2} - \mathrm{tr}'\,y) \\
    & \hspace{1em} + m \mathrm{tr}'\,a_{2l'-2} + \sum (\mathrm{tr}'\,x)(\mathrm{tr}'\,a_{2l'-2}) + \sum (\mathrm{tr}'\,x)(\mathrm{tr}'\,a_{2l'-3}) \\
    & \equiv 0 \pmod{J^{l'}}
  \end{split}\]
for some $m \in \Z$. (Remark that $\mathrm{tr}'\,a_{2l'-2} = \mathrm{tr}'\,a_{2l} \in J^{2l-1} \subset J^{l'}$ from the argument above.)
Therefore we see that Part (2) holds for $l'=l+1$.
This completes the proof of Proposition {\rmfamily \ref{P-7}}. $\square$

\begin{cor}\label{C-P-7}
With the notation above, we have
\begin{enumerate}
\item For any $k \geq 2$, $\mathrm{tr}'\,a_k^{-1} \in J^{k-1}$.
\item For any $l \geq 2$ and $b \in G$,
\[ \mathrm{tr}'\,ba_{2l-1}^{-1} - \mathrm{tr}'\,b, \,\,\, \mathrm{tr}'\,ba_{2l}^{-1} - \mathrm{tr}'\,b \in J^l. \]
\end{enumerate}
\end{cor}
\textit{Proof.}
This corollary is immediately proved by Proposition {\rmfamily \ref{P-7}} and
\[ [\alpha, \beta]^{-1} = [\beta \alpha \beta^{-1}, \beta^{-1}] \]
for any $\alpha, \beta \in G$.
This completes the proof of Corollary {\rmfamily \ref{C-P-7}}. $\square$

\vspace{0.5em}

\subsection{The structures of $\mathrm{gr}^k(J)$ for $G=F_n$, and $k=1$ and $2$}\label{Ss-Str}
\hspace*{\fill}\ 

\vspace{0.5em}

In this subsection, we assume that $G$ is a free group $F_n$ of rank $n \geq 2$. The goal of this subsection is to give a basis of $\mathrm{gr}^k(J)$ for $k=1$ and $2$.

\subsubsection{{\bf A basis of $\mathrm{gr}^1(J)$}}
\hspace*{\fill}\ 

\vspace{0.5em}

Here we show that the image of
\[\begin{split}
   T := \{ t_i' \,|\, & 1 \leq i \leq n \} \cup \{ t_{ij}' \,|\, 1 \leq i < j \leq n \}
   \cup \{ t_{ijk}' \,|\, 1 \leq i < j < k \leq n \} \subset \Q[t]
  \end{split}\]
by $\pi_{\Q}$ forms a basis of $\mathrm{gr}^1(J)$ as a $\Q$-vector space. To do this, it suffices to show
\begin{pro}\label{P-1}
For an ideal $J_0 := (T)$ of $\Q[t]$ generated by $T$, we have $I_{\Q} \subset J_0^2$.
\end{pro}
\textit{Proof.}
For any $f \in I_{\Q}$, set
\[\begin{split}
   f &:= \sum_{1 \leq i \leq n} a_i t_i' + \sum_{1 \leq i<j \leq n} a_{ij}t_{ij}' + \sum_{i<j<k} a_{ijk} t_{ijk}' \\
     & \hspace{3em} + (\text{terms of degree $\geq 2$ })
  \end{split}\]
for $a_i, a_{ij}, a_{ijk} \in \Q$.
Then, it suffices to show that $a_i=a_{ij}=a_{ijk}=0$.

\vspace{0.5em}

Choose any $i< j <k$ and fix it. Consider the interior
\[ D := \{ z \in \C \,|\, z \overline{z} < 1 \} \]
of the unit disk in $\C$. For any $s \in D$ and $l, m, t, u \in \C$, 
define a representation $\rho_1 : F_n \rightarrow \mathrm{SL}(2,\C)$ by
{\small
\[\begin{split}
   \rho_1(x_{i}) & := \begin{pmatrix} 1-s & 0 \\ 0 & (1-s)^{-1} \end{pmatrix},
        \,\,\, \rho_1(x_{j}) := \begin{pmatrix} 1-lt & l^2t \\ -t & 1+lt \end{pmatrix}, \,\,\,
   \rho_1(x_k) := \begin{pmatrix} 1-mu & m^2 u \\ -u & 1+mu \end{pmatrix}.
  \end{split}\]
}
If we consdier the power series expansion
\[ \frac{1}{1-s} = 1+s + s^2+ s^3+ \cdots \]
at the origin on $D$,
we can write $\mathrm{tr}'\, \rho_1(x_{i_1} \cdots x_{i_l})$ as a convergent power series of $s, t, u$:
\[\begin{split}
   \mathrm{tr}'\, \rho_1(x_i) & = \frac{s^2}{1-s} = s^2 + s^3 + s^4 + \cdots,  \\
   \mathrm{tr}'\, \rho_1(x_i x_j) &= \frac{1}{1-s} (s^2 + 2l st) \\
          & = s^2 + 2l st + s^3 + 2l s^2t + (\text{terms of degree $\geq 4$ }), \\
   \mathrm{tr}'\, \rho_1(x_j x_k) &= -(l-m)^2 tu, \\
   \mathrm{tr}'\, \rho_1(x_i x_j x_k) & = \frac{1}{1-s} \{ s^2 + 2lst + 2m su - (m-l)^2 tu + 2l(l-m)stu \\
     & \hspace{5em} -ls^2t -m s^2u +l(m-l)s^2tu \} \\
    & = s^2 + s^3  + 2lst + 2m su - (m-l)^2 tu + (l^2-m^2)stu +ls^2t + m s^2u \\
    & \hspace{1em}+ (\text{terms of degree $\geq 4$ }),
  \end{split}\]
and so on. This shows that $\mathrm{tr}'\, \rho_1(x_{i_1} \cdots x_{i_l})$ is eqaul to zero, or
the degrees of its monomials are greater than one.
Then we have
{\small
\begin{equation} \label{eq-marine} \begin{split}
   f(\mathrm{tr}'\,& \rho_1(x_{i_1} \cdots x_{i_l})) \\
     & = a_i (\mathrm{tr}'\, \rho_1(x_i)) + a_{jk} (\mathrm{tr}'\, \rho_1(x_j x_k))
         + \sum_{r<i} a_{ri} (\mathrm{tr}'\, \rho_1(x_r x_i)) + \sum_{i<r} a_{ir} (\mathrm{tr}'\, \rho_1(x_i x_r))  \\
     & \hspace{3em} + \sum_{j<k<p} a_{jkp} (\mathrm{tr}'\, \rho_1(x_j x_k x_p)) + \sum_{j<p<k} a_{jpk} (\mathrm{tr}'\, \rho_1(x_j x_p x_k)) \\
     & \hspace{3em} + \sum_{i \neq p<j<k} a_{pjk} (\mathrm{tr}'\, \rho_1(x_p x_j x_k)) \\
     & \hspace{3em} + \sum_{p<q<i} a_{pqi} (\mathrm{tr}'\, \rho_1(x_p x_q x_i)) +\sum_{p<i<q} a_{piq} (\mathrm{tr}'\, \rho_1(x_p x_i x_q)) \\
     & \hspace{3em} + \sum_{i<p<q} a_{ipq} (\mathrm{tr}'\, \rho_1(x_i x_p x_q)) \\
     & =0.
  \end{split}\end{equation}}
By the uniqueness of the power series expansion on $D$,
each of the coefficients of the monomials in $f(\mathrm{tr}'\, \rho_1(x_{i_1} \cdots x_{i_l}))$ must be equal to zero.
Here we observe the coefficients of the monomials of degree less than four.

\vspace{0.5em}

First, from the coefficient of $stu$, we obtain $(l^2-m^2) a_{ijk}=0$. Since we can choose $l, m \in \C$ arbitrary, we see $a_{ijk}=0$.
Therefore (\ref{eq-marine}) reduces
\[\begin{split}
   f(\mathrm{tr}'\, \rho_1(x_{i_1} \cdots x_{i_l}))
     & = a_i (\mathrm{tr}'\, \rho_1(x_i)) + a_{jk} (\mathrm{tr}'\, \rho_1(x_j x_k)) \\
     & \hspace{3em} + \sum_{r<i} a_{ri} (\mathrm{tr}'\, \rho_1(x_r x_i)) + \sum_{i<r} a_{ir} (\mathrm{tr}'\, \rho_1(x_i x_r)) \\
     & = 0.
  \end{split}\]
Next, from the coefficient of $st$, $su$ and $tu$, we see $a_{ij}=0$, $a_{ik}=0$ and $a_{jk}=0$ respectively. Hence we have
\[ f(\mathrm{tr}'\, \rho_1(x_{i_1} \cdots x_{i_l})) = a_i (\mathrm{tr}'\, \rho_1(x_i)) =0. \]
Furthermore, from the coefficient of $s^2$, we see $a_{i}=0$.

\vspace{0.5em}

On the other hand, for any $t, u \in D$,
define a representation $\rho_1', \rho_1'' : F_n \rightarrow \mathrm{SL}(2,\C)$ by
{\small
\[\begin{split}
   \rho_1'(x_{r}) & :=\begin{cases}
                        \begin{pmatrix} 1-t & 1 \\ 0 & (1-t)^{-1} \end{pmatrix} \hspace{1em} & \text{if} \hspace{1em} r=j, \\
                        E_2 & \text{if} \hspace{1em} r \neq j,
                     \end{cases}\\
   \rho_1''(x_{r}) & :=\begin{cases}
                        \begin{pmatrix} 1-u & 1 \\ 0 & (1-u)^{-1} \end{pmatrix} \hspace{1em} & \text{if} \hspace{1em} r=k, \\
                        E_2 & \text{if} \hspace{1em} r \neq k.
                     \end{cases}
  \end{split}\]
}
Then, by an argument similar to that in the above, from the coefficients of $t^2$ in $f(\mathrm{tr}'\, \rho_1'(x_{i_1} \cdots x_{i_l}))$,
and of $u^2$ in $f(\mathrm{tr}'\, \rho_1''(x_{i_1} \cdots x_{i_l}))$, we obtain $a_j=0$ and $a_k=0$ respectively.

\vspace{0.5em}

Therefore we conclude that $f \in J_0^2$. This completes the proof of Proposition {\rmfamily \ref{P-1}}. $\square$

\subsubsection{{\bf A basis of $\mathrm{gr}^2(J)$}}
\hspace*{\fill}\ 

\vspace{0.5em}

Set
\[\begin{split}
   S_1 := & \{ t_i' t_j' \,|\, 1 \leq i \leq j \leq n \} \cup \{ t_i' t_{ab}' \,|\, 1 \leq i \leq n, \,\, 1 \leq a < b \leq n \} \\
        & \cup \{ t_i' t_{abc}' \,|\, 1 \leq i \leq n, \,\, 1 \leq a<b<c \leq n \} \\
        & \cup \{ t_{ij}' t_{ab}' \,|\, 1 \leq i<j \leq n, \,\, 1 \leq a<b \leq n, \,\, (i,j) \leq (a,b) \}, \\
   S_2 := & \{ t_{ab}' t_{abc}', \, t_{ac}' t_{abc}', \, t_{bc}' t_{abc}' \,|\, 1 \leq a<b<c \leq n \} \\
        & \cup \{ t_{ia}' t_{abc}', t_{ib}' t_{abc}', t_{ic}' t_{abc}', t_{ia}' t_{ibc}', t_{ab}' t_{iac}', t_{ab}' t_{ibc}', t_{ac}' t_{ibc}', t_{ib}' t_{iac}' \,|\, 1 \leq i < a< b<c \leq n \} \\
        & \cup \{ t_{ja}' t_{ibc}', t_{jb}' t_{iac}', t_{jc}' t_{iab}', t_{ab}' t_{ijc}', t_{ac}' t_{ijb}', t_{bc}' t_{ija}'
            \,|\, 1 \leq i<j<a<b<c \leq n \}
  \end{split}\]
and $S := S_1 \cup S_2$.
Here we show that $\pi_{\Q}(S)$ forms a basis of $\mathrm{gr}^2(J)$ as a $\Q$-vector space.

\vspace{0.5em}

First, we show
\begin{pro}\label{P-haruhi}
$\pi_{\Q}(S)$ generates $\mathrm{gr}^2(J)$. 
\end{pro}
\textit{Proof.}
Set
\[\begin{split}
   S' & := \{ t_{ij}' t_{abc}' \,|\, 1 \leq i < j \leq n, \,\,\, 1 \leq a<b<c \leq n \}, \\
   S'' & := \{ t_{ijk}' t_{abc}' \,|\, 1 \leq i<j<k \leq n, \,\, 1 \leq a<b<c \leq n, \,\,\, (i,j,k) \leq (a,b,c) \}. 
  \end{split}\]
Then $\mathrm{gr}^2(J)$ is generated by $\pi_{\Q}(S_1 \cup S' \cup S'')$.
Consider relations
\begin{eqnarray}
 & & (2t_{ijk} -t_i t_{jk} - t_j t_{ik} - t_k t_{ij} + t_i t_j t_k) (2t_{abc} -t_a t_{bc} - t_b t_{ac} - t_c t_{ab} + t_a t_b t_c) \\
 & & \hspace{5em} = \begin{vmatrix} t_i & t_{ia} & t_{ib} & t_{ic} \\ t_j & t_{ja} & t_{jb} & t_{jc} \\
                                    t_k & t_{ka} & t_{kb} & t_{kc} \\ 2 & t_a & t_b & t_c \end{vmatrix} \nonumber
\end{eqnarray}
in $\Q[t]/I_{\Q}$ where $t_{ii}=t_i^2 -2$. (For details, see Corollary 4.12 in \cite{AJM}.)
Substituting $t_{i_1 \cdots i_l}' +2$ to each of $t_{i_1 \cdots i_l}$ in the equations above, we verify that
$\pi_{\Q}(t_{ijk}' t_{abc}')$ is written as a polynomial of the indeterminates $t_i'$ and $t_{ij}'$ in $\Q[t]/I_{\Q}$.
Hence we see that $\pi_{\Q}(S\cup S')$ generates $\mathrm{gr}^2(J)$.

\vspace{0.5em}

Next, we reduce the generators of $\pi_{\Q}(S')$. Consider a quotient $\Q$-vector space $V$ of $\mathrm{gr}^2(J)$ by a subspace generated by
$\pi_1(S_1)$. We write the equality in $V$ as $\doteq$.
Now, fix $1 \leq i < a< b< c \leq n$. Then we have elements
\[\begin{split}
   p_2 &:= t_i t_{abc} +t_{acb}-t_{abc} - t_{ia}t_{ibc} + t_{ib}t_{iac} - t_{ic} t_{iab} - t_{i} t_b t_{iac} + t_b t_{ia} t_{ic}, \\
   p_3 &:= t_{ib} t_{iabc} - t_{iab} t_{ibc} - t_i t_{iac} + t_{ia} t_{ic} - t_a t_c + 2t_{ac} - t_b t_{abc} + t_{ab}t_{bc}, \\
   p_4 &:= t_{iba} t_{iabc} - t_{ia}t_{ab}t_{ibc} - t_{ab}t_{abc} + t_it_b t_{ibc} - t_{ib}t_{ibc} - t_{ia}t_{iac}
           + t_a t_{ab} t_{bc} + t_a t_{ia} t_{ic} \\
       & \hspace{4em} + t_a t_{ac} - t_i t_{ic} - t_b t_{bc} + 2t_c, \\
   (p_3)^{\sigma_{bc}} &:= t_{ic} t_{iacb} -t_{iac}t_{icb} - t_it_{iab} + t_{ia}t_{ib} - t_at_b +2 t_{ab} -t_c t_{acb} + t_{ac}t_{cb}
  \end{split}\]
in $I$ due to Whittemore \cite{Whi}, where $\sigma_{bc} \in \mathrm{Aut}\,F_n$ is an automorphism such that
\[ x_r \mapsto \begin{cases}
                  x_c \hspace{1em} & \mathrm{if} \hspace{1em} r=b, \\
                  x_b & \mathrm{if} \hspace{1em} r=c, \\
                  x_r & \mathrm{if} \hspace{1em} r \neq b, c.
               \end{cases}\]
From the above, by using (\ref{eq-20}) and a straightforward calculation, we obtain equations
\begin{eqnarray}
  t_{ic}' t_{iab}' & \doteq & -t_{ia}' t_{ibc}' + t_{ib}' t_{iac}', \\
  t_{bc}' t_{iab}' & \doteq & - t_{ab}' t_{ibc}' + t_{ib}' t_{abc}' - (t_{ic}' t_{iab}' + t_{ia}' t_{ibc}' - t_{ib}' t_{iac}') \nonumber \\
   & \doteq & - t_{ab}' t_{ibc}' + t_{ib}' t_{abc}', \\
  t_{ac}' t_{iab}' & \doteq & t_{ia}' t_{abc}' + t_{ab}' t_{iac}' - (t_{ic}' t_{iab}' + t_{ia}' t_{ibc}' - t_{ib}' t_{iac}') \nonumber \\
   & & \hspace{2em}   +(- t_{bc} t_{iab}' - t_{ab}' t_{ibc}' + t_{ib}' t_{abc}') \nonumber \\
   & \doteq &  t_{ia}' t_{abc}' + t_{ab}' t_{iac}', \\
  t_{bc}' t_{iac}' & \doteq & -t_{ic}' t_{abc}' + t_{ac}' t_{ibc}' + (t_{ic}' t_{iab}' + t_{ia}' t_{ibc}' - t_{ib}' t_{iac}') \nonumber \\
     & \doteq & -t_{ic}' t_{abc}' + t_{ac}' t_{ibc}'
\end{eqnarray}
in $V$ respectively. Hence we can remove $t_{ic}' t_{iab}'$, $t_{bc}' t_{iab}'$, $t_{ac}' t_{iab}'$ and $t_{bc}' t_{iac}'$ from
the generating set $\pi_{\Q}(S')$.

\vspace{0.5em}

Fix $1 \leq i < j< a< b< c \leq n$. Using (\ref{eq-5}), we have
\[\begin{split} 2t_{(ij)abc} & = t_{ij} t_{abc} + t_a t_{ijbc} + t_b t_{ijac} + t_c t_{ijab} + t_{ija} t_{bc}
                   - t_{ijb} t_{ac} + t_{ijc} t_{ab} \\
    & \hspace{1em} - t_{ij} t_a t_{bc} - t_a t_b t_{ijc} - t_{ij} t_c t_{ab} -t_b t_c t_{ija} + t_{ij} t_a t_b t_c
  \end{split}\]
in $\Z[t]$. On the other hand, we have
\[\begin{split} 2t_{(ja)bci} & = t_{ja} t_{bci} + t_b t_{jaci} + t_c t_{jabi} + t_i t_{jabc} + t_{jab} t_{ci}
                   - t_{jac} t_{bi} + t_{jai} t_{bc} \\
    & \hspace{1em} - t_{ja} t_b t_{ci} - t_b t_c t_{jai} - t_{ja} t_i t_{bc} -t_c t_i t_{jab} + t_{ja} t_b t_c t_i.
  \end{split}\]
Hence, from $2t_{(ij)abc}=2t_{(ja)bci}$, by using (\ref{eq-5}) again, we obtain
\begin{equation}\label{eq-40}
 t_{ij}' t_{abc}' + t_{ib}' t_{jac}' - t_{ic}' t_{jab}' \doteq t_{ja}' t_{ibc}' - t_{ab}' t_{ijc}' + t_{ac}' t_{ijb}'.
\end{equation}
Similarly, from equations $2t_{(ij)abc}=2t_{(ab)cij}$, $2t_{(ij)abc}=2t_{(bc)ija}$ and $2t_{(ij)abc}=2t_{(ci)jab}$,
we obtain
\begin{eqnarray}
 t_{ic}' t_{jab}' & \doteq & t_{jc}' t_{iab}' - t_{ac}' t_{ijb}' + t_{bc}' t_{ija}', \label{eq-41} \\
 t_{ia}' t_{jbc}' & \doteq & t_{ja}' t_{ibc}' - t_{ab}' t_{ijc}' + t_{ac}' t_{ijb}', \label{eq-42} \\
 t_{ij}' t_{abc}' - t_{ic}' t_{jab}' & \doteq & t_{ja}' t_{ibc}' - t_{jb}' t_{iac}' + t_{ac}' t_{ijb}' - t_{bc}' t_{ija}' \label{eq-43}
\end{eqnarray}
respectively. From (\ref{eq-40}), (\ref{eq-41}) and (\ref{eq-43}), we see
\begin{eqnarray}
 t_{ib}' t_{jac}' & \doteq & - t_{ab}' t_{ijc}' + t_{jb}' t_{iac}' + t_{bc}' t_{ija}' \label{eq-44} \\
 t_{ij}' t_{abc}' & \doteq & t_{ja}' t_{ibc}' - t_{jb}' t_{iac}' + t_{jc}' t_{iab}'. \label{eq-45}
\end{eqnarray}
Therefore, by observing (\ref{eq-45}), (\ref{eq-42}), (\ref{eq-44}) and (\ref{eq-41}), we can remove
$t_{ij}' t_{abc}'$, $t_{ia}' t_{jbc}'$, $t_{ib}' t_{jac}'$ and $t_{ic}' t_{jab}'$ from
the generating set $\pi_{\Q}(S')$.

\vspace{0.5em}

Then we obtain the required result. This completes the proof of Proposition {\rmfamily \ref{P-haruhi}}. $\square$

\vspace{0.5em}

Next we prove
\begin{pro}\label{P-2}
Elements in $\pi_{\Q}(S)$ are linearly independent in $\mathrm{gr}^2(J)$.
\end{pro}
\textit{Proof.}
Set
\[\begin{split}
   g &:= \sum_{1 \leq i \leq j \leq n} d_{i,j} t_i' t_j' + \sum_{\substack{1 \leq i \leq n, \\[1pt] 1 \leq a < b \leq n}} d_{i, ab} t_{i}'t_{ab}'
         + \sum_{\substack{1 \leq i \leq n, \\[1pt] 1 \leq a < b < c \leq n}} d_{i, abc} t_{i}'t_{abc}' \\
     & \hspace{2em} + \sum_{(1,2) \leq (i,j) \leq (a,b) \leq (n-1,n) } d_{ij, ab} t_{ij}' t_{ab}'
             + \sum_{\substack{1 \leq i < j \leq n, \,\, 1 \leq a < b < c \leq n; \\[1pt] t_{ij}' t_{abc}' \in S_2 }} d_{ij, abc} t_{ij}' t_{abc}' \in \Q[t]
  \end{split}\]
for $d_{i,j}, d_{i,ab}, d_{i,abc}, d_{ij,ab}, d_{ij,abc} \in \Q$. Assume $\pi_{\Q}(g) \in J^3$.
Then, it suffices to show that $d_{i,j}=d_{i,ab}=d_{i,abc}=d_{ij,ab}=d_{ij,abc}=0$.

\vspace{0.5em}

{\bf Step 1.} $d_{ij, abc}=0$ for any $1 \leq i<j \leq n$ and $1 \leq a<b<c \leq n$ such that $t_{ij}' t_{abc}' \in S_2$.

\vspace{0.5em}

Set $N_1:= \{ i, j, a, b, c \}$. We consider three cases according to the number of elements in $N_1$.

\vspace{0.5em}

{\bf Case 1-1.} $\sharp N_1 =3$.

\vspace{0.5em}

Assume $N_1=\{ a, b, c \}$ and $a<b<c$. We show
\begin{equation}
    d_{ab, abc} =d_{ac, abc}= d_{bc, abc} =0.
\label{case1-1} \end{equation}

To do this, for any $k, l, m, s, t, u \in \C$, 
consider a representation $\rho_2 : F_n \rightarrow \mathrm{SL}(2,\C)$ defined by
{\small
\[\begin{split}
   \rho_2(x_{a}) & := \begin{pmatrix} 1-ks & k^2 s \\ -s & 1+ks \end{pmatrix}, \,\,\, \rho_2(x_{b}) := \begin{pmatrix} 1-lt & l^2t \\ -t & 1+lt \end{pmatrix}, \,\,\,
      \rho_2(x_c) := \begin{pmatrix} 1-mu & m^2 u \\ -u & 1+mu \end{pmatrix}
  \end{split}\]
}
and $\rho_2(x_r)= E_2$ for $r \neq a, b, c$. Then we have
\[\begin{split}
   \mathrm{tr}'\, \rho_2(x_a x_b) &= -(k-l)^2 st, \\ \mathrm{tr}'\, \rho_2(x_a x_b x_c) & = (k-l)(l-m)(m-k) stu -(k-l)^2st - (l-m)^2 tu -(m-k)^2su,
  \end{split}\]
and so on.

\vspace{0.5em}

Consider each of $\mathrm{tr}'\, \rho_2(x_{i_1} \cdots x_{i_l})$ as a polynomial of $s, t, u$ with rational coefficients.
Then by the observation above, each of them is zero, or is of degree greater than one.
Hence, since $\pi_{\Q}(g) \in J^3$, the degree of $g(\mathrm{tr}'\, \rho_2(x_{i_1} \cdots x_{i_l}))$ must be greater than five.
This shows that each of the coefficients of
$s^2 t^2 u, s^2 tu^2, st^2u^2$ in $g(\mathrm{tr}'\, \rho_2(x_{i_1} \cdots x_{i_l}))$ is equal to zero. Hence we see
\[\begin{split}
    & -(k-l)^2 (k-l)(l-m)(m-k) d_{ab, abc} -(k-m)^2 (k-l)(l-m)(m-k) d_{ac, abc} \\
    & - (l-m)^2 (k-l)(l-m)(m-k) d_{bc, abc} \\
    & =0.
  \end{split}\]
Furthermore, since we can choose $k, l, m \in \C$ arbitrary, we obtain (\ref{case1-1}).

\vspace{0.5em}

{\bf Case 1-2.} $\sharp N_1 =4$.

\vspace{0.5em}

Assume $N_1=\{ i, a, b, c \}$ and $i<a<b<c$. It suffices to show
\begin{equation}\begin{split}
   d_{ia, abc} =d_{ab, iac} =0, \hspace{1em} & d_{ib, abc} =d_{ab, ibc} =0, \\
   d_{ic, abc} =d_{ac, ibc} =0, \hspace{1em} & d_{ia, ibc} =d_{ib, iac} =0. 
\end{split} \label{case1-3} \end{equation}

To do this, for any $k, l, m, p, s, t, u, v \in \C$, 
consider a representation $\rho_3 : F_n \rightarrow \mathrm{SL}(2,\C)$ defined by

\[\begin{split}
   \rho_3(x_{i}) := \begin{pmatrix} 1-pv & p^2 v \\ -v & 1+pv \end{pmatrix}
  \end{split}\]
and $\rho_3(x_r):= \rho_2(x_r)$ for $r \neq i$. 
From the coefficient of $s^2 tuv$ in $g(\mathrm{tr}'\, \rho_3(x_{i_1} \cdots x_{i_l}))$, we have
\[\begin{split}
    -(k-p)^2 (k-l)(l-m)(m-k) d_{ia, abc} - (k-l)^2 (p-k)(k-m)(m-p) d_{ab, iac} = 0
  \end{split}\]
and from the coefficients of $k^4l, k^4p$, we see
\[ d_{ia, abc}=d_{ab, iac}=0. \]

\vspace{0.5em}

Similarly, by observing the coefficients of $st^2uv$, $stu^2v$ and $stuv^2$ in $g(\mathrm{tr}'\, \rho_3(x_{i_1} \cdots x_{i_l}))$, we see
\[\begin{split}
 -(l-p)^2 & (k-l)(l-m)(m-k) d_{ib, abc} - (l-k)^2 (p-l)(l-m)(m-p) d_{ab, ibc} =0, \\
 -(p-m)^2 & (k-l)(l-m)(m-k) d_{ic, abc} - (k-m)^2 (p-l)(l-m)(m-p) d_{ac, ibc} =0, \\
 -(p-k)^2 & (p-l)(l-m)(m-p) d_{ia, ibc} - (p-l)^2 (p-k)(k-m)(m-p) d_{ib, iac} =0
\end{split}\]
respectively. From these, we obtain (\ref{case1-3}).

\vspace{0.5em}

{\bf Case 1-3.} $\sharp N_1 =5$.

\vspace{0.5em}

Assume $i<j<a<b<c$. Then it suffices to show
\begin{equation}\begin{split}
    d_{ja, ibc} & = d_{jb, iac}=d_{jc, iab} = d_{ab, ijc} = d_{ac, ijb} = d_{bc, ija} =0.
  \end{split} \label{case1-2} \end{equation}

To begin with, for any $k, l, m, p, q, s, t, u, v, w \in \C$, 
consider a representation $\rho_4 : F_n \rightarrow \mathrm{SL}(2,\C)$ defined by
\[\begin{split}
   \rho_4(x_{j}) := \begin{pmatrix} 1-qw & q^2 w \\ -w & 1+qw \end{pmatrix}
  \end{split}\]
and $\rho_4(x_r):= \rho_3(x_r)$ for $r \neq j$. 

\vspace{0.5em}

Then $g(\mathrm{tr}'\, \rho_4(x_{i_1} \cdots x_{i_l})) \in \C$ is written as a polynomial of $k, l, m, p, q$, $s, t, u, v, w$ with rational coefficients.
By observing the coefficient of $stuvw$ in $g(\mathrm{tr}'\, \rho_4(x_{i_1} \cdots x_{i_l}))$, we have
\[\begin{split}
    & - (k-q)^2 (p-l)(l-m)(m-p) d_{ja, ibc} - (l-q)^2 (p-k)(k-m)(m-p) d_{jb, iac} \\
    & - (q-m)^2 (p-k)(k-l)(l-p) d_{jc, iab} - (k-l)^2 (p-q)(q-m)(m-p) d_{ab, ijc} \\
    & - (k-m)^2 (p-q)(q-l)(l-p) d_{ac, ijb} - (l-m)^2 (p-q)(q-k)(k-p) d_{bc, ija} \\
    & =0.
  \end{split}\]

Furthermore, from the coefficients of $p^2 m^2 k, p^2 m^2 l$, $p^2l^2m, p^2q^2l$, $p^2q^2m, p^2l^2q, p^2k^2m$, we obtain
\[\begin{split}
 d_{bc, ija} &= d_{jc, iab}, \hspace{1em} d_{ac, ijb} = - d_{jc, iab}, \hspace{1em} d_{ab, ijc} = - d_{jb, iac}, \\
 d_{jc, iab} &= d_{ja, ibc}, \hspace{1em} d_{jb, iac} = - d_{ja, ibc}, \hspace{1em} d_{ab, ijc} = - d_{bc, ija}, \hspace{1em} d_{ja, ibc} =- d_{ab, ijc}
\end{split}\]
respectively. From this, we obtain (\ref{case1-2}).

\vspace{0.5em}

{\bf Step 2.} $d_{ij, ab}=0$ for any $1 \leq i<j \leq n$ and $1 \leq a<b \leq n$ and $(i,j) \leq (a,b)$.

\vspace{0.5em}

Set $N_2:= \{ i, j, a, b \}$. We consider three cases according to the number of elements in $N_2$.

\vspace{0.5em}

First, we consider the case where $\sharp N_2=2$. We show $d_{ab, ab}=0$ for any $1 \leq a < b \leq n$.
Recall the representation $\rho_2$. By observing the coefficients of $s^2t^2$ in $g(\mathrm{tr}'\, \rho_2(x_{i_1} \cdots x_{i_l}))$, we obtain
\[ (k-l)^4 d_{ab, ab} = 0. \]
This shows that $d_{ab, ab}=0$.

\vspace{0.5em}

Next, consider the case where $\sharp N_2=3$. It suffices to show that for any $1 \leq i< a < b \leq n$,
\begin{equation}
 d_{ia, ab}= d_{ib, ab} = d_{ia, ib} = 0. \label{haruhi-12}
\end{equation}
By observing the coefficients of $s^2tv$, $st^2v$ and $stv^2$ in $g(\mathrm{tr}'\, \rho_3(x_{i_1} \cdots x_{i_l}))$, we obtain
\[\begin{split}
   (p-k)^2(k-l)^2 d_{ia, ab} & = 0, \hspace{1em} (p-l)^2(k-l)^2 d_{ib, ab} =0, \\
   (p-k)^2(p-l)^2 d_{ia, ib} & =0. 
  \end{split}\]
From the coefficients of $k^4$, $l^4$ and $p^4$, we see (\ref{haruhi-12}).

\vspace{0.5em}

Finally, consider the case where $\sharp N_2=4$. It suffices to show that for any $1 \leq i< j< a < b \leq n$,
\begin{equation}
 d_{ij, ab}= d_{ia, jb} = d_{ib, ja} = 0. \label{haruhi-13}
\end{equation}
By observing the coefficients of $stuv$ in $g(\mathrm{tr}'\, \rho_4(x_{i_1} \cdots x_{i_l}))$, we obtain
\[\begin{split}
   (p-q)^2(k-l)^2 d_{ij, ab} + (p-k)^2(q-l)^2 d_{ia, jb} =0 + (p-l)^2(q-k)^2 d_{ib, ja} =0. 
  \end{split}\]
From the coefficients of $pqk^2$, $qlk^2$ and $plq^2$, we see (\ref{haruhi-13}).

\vspace{0.5em}

{\bf Step 3.} $d_{i, abc}= 0$ for any $1 \leq i \leq n$ and $1 \leq a<b<c \leq n$.

\vspace{0.5em}

First, assume $i \neq a, b, c$.
For any $v \in D$ and $k, l, m, s, t, u \in \C$, 
consider a representation $\rho_5 : F_n \rightarrow \mathrm{SL}(2,\C)$ defined by
\[\begin{split}
   \rho_5(x_{i}) := \begin{pmatrix} 1-v & 0 \\ 0 & (1-v)^{-1} \end{pmatrix}
  \end{split}\]
and $\rho_5(x_r):= \rho_3(x_r)$ for $r \neq i$. 
By observing the coefficients of $stu v^2$ in $g(\mathrm{tr}'\, \rho_5(x_{i_1} \cdots x_{i_l}))$, we obtain
\[ (k-l)(l-m)(m-k) d_{i, abc} = 0. \]
This shows that $d_{i, abc}=0$.

\vspace{0.5em}

Next, we consider the case where $i= a, b$ or $c$. For any $s \in D$ and $l, m, t, u \in \C$, 
consider a representation $\rho_6 : F_n \rightarrow \mathrm{SL}(2,\C)$ defined by
\[\begin{split}
   \rho_6(x_a) := \begin{pmatrix} 1-s & 0 \\ 0 & (1-s)^{-1} \end{pmatrix}
  \end{split}\]
and $\rho_6(x_r):= \rho_2(x_r)$ for $r \neq a$. Then from the coefficients of $s^3tu$ in $g(\mathrm{tr}'\, \rho_6(x_{i_1} \cdots x_{i_l}))$,
we have
\[ (l^2-m^2) d_{a, abc} -(l-m)^2 d_{a, bc} = 0. \]
This shows $d_{a, abc} = d_{a, bc} =0$. Similarly, we can obtain
\[ d_{b, abc} = d_{c, abc} =0. \]

\vspace{0.5em}

{\bf Step 4.} $d_{i, ab}=0$ for any $1 \leq i \leq n$ and $1 \leq a<b \leq n$.

\vspace{0.5em}

Assume $i \neq a, b$.
By observing the coefficients of $st v^2$ in $g(\mathrm{tr}'\, \rho_5(x_{i_1} \cdots x_{i_l}))$, we obtain
\[ -(k-l)^2 d_{i, ab} = 0. \]
This shows $d_{i, ab}=0$. (Remark that for $i<a<b$, this has already been obtained in Step 3.)

\vspace{0.5em}

On the other hand, from the coefficients of $s^3t$ in $g(\mathrm{tr}'\, \rho_6(x_{i_1} \cdots x_{i_l}))$,
we have
\[ 2l d_{a, ab} = 0. \]
This shows $d_{a, ab} = 0$. Similarly, we can obtain
\[ d_{b, ab} =0. \]

\vspace{0.5em}

{\bf Step 5.} $d_{i,a}=0$ for any $1 \leq i, a \leq n$.

\vspace{0.5em}

Assume $i \neq a$.
For any $s, v \in D$, 
consider a representation $\rho_7 : F_n \rightarrow \mathrm{SL}(2,\C)$ defined by
\[\begin{split}
   \rho_7(x_i) := \begin{pmatrix} 1-s & 0 \\ 0 & (1-s)^{-1} \end{pmatrix}, \,\,\, \rho_7(x_a) := \begin{pmatrix} 1-v & 0 \\ 0 & (1-v)^{-1} \end{pmatrix}
  \end{split}\]
and $\rho_7(x_r):= E_2$ for $r \neq i, a$. 
Then by observing the coefficients of $s^4$, $s^2v^2$ and $v^4$ in $g(\mathrm{tr}'\, \rho_{7}(x_{i_1} \cdots x_{i_l}))$, we obtain
\[ d_{i,i} = 0, \,\,\, d_{i,a}=0, \,\,\, \text{and} \,\,\, d_{a,a}=0 \]
respectively.

\vspace{0.5em}

Therefore we have obtained all coefficients of $g$ are eqaul to zero.
This compeletes the proof of Proposition {\rmfamily \ref{P-2}}. $\square$

\section{A central filtration $\mathcal{E}_G(k)$}\label{S-AutG}

\vspace{0.5em}

In this section, for any group $G$, we introduce a descending filtration of $\mathrm{Aut}\,G$ consisting of its normal subgroups.
This is an analogue of the Andreadakis-Johnson filtration of $\mathrm{Aut}\,G$. (For details for the Andreadakis-Johnson filtration,
see \cite{S06} or \cite{S14}, for example.)

\subsection{Definition of $\mathcal{E}_G(k)$}
\hspace*{\fill}\ 

\vspace{0.5em}

For any $k \geq 1$, let
\[ \mathcal{E}_G(k) := \mathrm{Ker}(\mathrm{Aut}\,G \rightarrow \mathrm{Aut}(J/J^{k+1})) \]
be the kernel of a homomorphism $\mathrm{Aut}\,G \rightarrow \mathrm{Aut}(J/J^{k+1})$ which is induced from the action of $\mathrm{Aut}\,G$
on $J/J^{k+1}$. Then the groups $\mathcal{E}_G(k)$ define a descending filtration
\[ \mathcal{E}_G(1) \supset \mathcal{E}_G(2) \supset \cdots \supset \mathcal{E}_G(k) \supset \cdots \]
of $\mathrm{Aut}\,G$. Here we show that this is a central filtration.

\vspace{0.5em}

For any $f \in J$ and $\sigma \in \mathrm{Aut}\,G$, set
\[ s_{\sigma}(f) := f^{\sigma}- f \in J. \]
Then we have
\begin{lem}\label{L-1}
For any $f \in J$ and $\sigma, \tau \in \mathrm{Aut}\,G$,
\begin{enumerate}
\item $s_{\sigma \tau}(f) = (s_{\sigma}(f))^{\tau} + s_{\tau}(f)$,
\item $s_{1_G}(f) = 0$,
\item $s_{\sigma^{-1}}(f) = - (s_{\sigma}(f))^{\sigma^{-1}}$,
\item $s_{[\sigma, \tau]}(f) = \{ s_{\tau}(s_{\sigma}(f)) - s_{\sigma}(s_{\tau}(f)) \}^{\sigma^{-1}\tau^{-1}}$.
\end{enumerate}
\end{lem}
\textit{Proof.}
The part of (1), (2) and (3) is straightforward. Here we prove the part (4). Using (1), (2) and (3), we obtain
\[\begin{split}
   s_{[\sigma, \tau]}(f) & \stackrel{(1)}{=} (s_{\sigma \tau}(f))^{\sigma^{-1} \tau^{-1}} + s_{\sigma^{-1} \tau^{-1}}(f) \\
    & \stackrel{(3)}{=} (s_{\sigma \tau}(f))^{\sigma^{-1} \tau^{-1}} - (s_{\tau \sigma}(f))^{\sigma^{-1} \tau^{-1}} \\
    & \stackrel{(1)}{=} \{ (s_{\sigma}(f))^{\tau} + s_{\tau}(f)  - (s_{\tau}(f))^{\sigma} - s_{\sigma}(f) \}^{\sigma^{-1} \tau^{-1}} \\
    & = \{ s_{\tau}(s_{\sigma}(f)) - s_{\sigma}(s_{\tau}(f)) \}^{\sigma^{-1}\tau^{-1}}.
  \end{split}\]
This completes the proof of Lemma {\rmfamily \ref{L-1}}. $\square$

\begin{lem}\label{L-2}
For any $k, l \geq 1$, $f \in J^l$ and $\sigma \in \mathcal{E}_G(k)$, we have $s_{\sigma}(f) \in J^{k+l}$.
\end{lem}
\textit{Proof.}
It suffices to show the lemma for the case where $f$ is (the coset class of) a monomial
$t_{a_1 \cdots a_{r_1}}' t_{b_1 \cdots b_{r_2}}' \cdots t_{c_1 \cdots c_{r_l}}'$. Then we have
\[\begin{split}
   s_{\sigma}(f) & = f^{\sigma} - f \\
       & = (t_{a_1 \cdots a_{r_1}}')^{\sigma} \cdots (t_{c_1 \cdots c_{r_l}}')^{\sigma} - t_{a_1 \cdots a_{r_1}}' \cdots t_{c_1 \cdots c_{r_l}}' \\
        & = (t_{a_1 \cdots a_{r_1}}' + s_{\sigma}(t_{a_1 \cdots a_{r_1}}')) \cdots (t_{c_1 \cdots c_{r_l}}' + s_{\sigma}(t_{c_1 \cdots c_{r_l}}'))
              - t_{a_1 \cdots a_{r_1}}' \cdots t_{c_1 \cdots c_{r_l}}'.
  \end{split}\]
By the definition of $\mathcal{E}_G(k)$, the elements $s_{\sigma}(t_{a_1 \cdots a_{r_1}}'), \ldots, s_{\sigma}(t_{c_1 \cdots c_{r_l}}')$ belong to
$J^{k+1}$. Therefore, we obtain $s_{\sigma}(f) \in J^{k+l}$. This completes the proof of Lemma {\rmfamily \ref{L-2}}. $\square$

\begin{pro}\label{P-3}
For any $k, l \geq 1$, $[\mathcal{E}_G(k), \mathcal{E}_G(l)] \subset \mathcal{E}_G(k+l)$.
\end{pro}
\textit{Proof.}
For any $\sigma \in \mathcal{E}_G(k)$, $\tau \in \mathcal{E}_G(l)$ and $f \in J$, by Lemmas {\rmfamily \ref{L-1}} and {\rmfamily \ref{L-2}}, we see
\[\begin{split}
   s_{[\sigma, \tau]}(f) & = \{ s_{\tau}(s_{\sigma}(f)) - s_{\sigma}(s_{\tau}(f)) \}^{\sigma^{-1}\tau^{-1}} \\
    & \equiv 0 \pmod{J^{k+l+1}}. \\
  \end{split}\]
Hence $[\sigma, \tau] \in \mathcal{E}_G(k+l)$. This completes the proof of Proposition {\rmfamily \ref{P-3}}. $\square$

\vspace{0.5em}

This proposition shows that the filtration $\mathcal{E}_G(k)$ is a central filtration of $\mathrm{Aut}\,G$.
Next, our interests is how different the filtration $\mathcal{E}_G(k)$ is from the Andreadakis-Johnson filtration $\mathcal{A}_G(k)$.
We consider this problem for the case where $G=F_n$.

\vspace{0.5em}

\subsection{The group $\mathcal{E}_{F_n}(1)$}
\hspace*{\fill}\ 

\vspace{0.5em}

Here we show that $\mathcal{E}_{F_n}(1) = \mathrm{Inn}\,F_n \cdot \mathcal{A}_{F_n}(2)$. First, we show that $\mathcal{E}_{F_n}(1)$ is contained in
the IA-automorphism group $\mathrm{IA}_n = \mathcal{A}_{F_n}(1)$. In the following, we always identify $t_{i_1 \cdots i_l}' \in \Q[t]/I_{\Q}$ with
$\mathrm{tr}'\,x_{i_1} \cdots x_{i_l} \in \mathcal{F}(R(G),\C)$ through $\pi_{\Q}$.

\vspace{0.5em}

To begin with, we prepare some lemmas.
\begin{lem}\label{L-3}
For any $s \in \C$, Set $A := \begin{pmatrix} s+2 & 1 \\ -1 & 0 \end{pmatrix}$, and
\[ \mathrm{tr}\,A^m = x_m^{(0)} + x_m^{(1)}s + x_m^{(2)}s^2 + \cdots \]
for any $m \in \Z$. Then $x_m^{(0)} =2$ and $x_m^{(1)}=m^2$.
\end{lem}
This lemma is obtained by a straightforward calculation.

\begin{lem}\label{L-4}
For any $1 \leq i \leq n$ and a word $w := x_{i_1}^{e_1} \cdots x_{i_l}^{e_l} \in F_n$, assume
\[ \mathrm{tr}'(w) \equiv \mathrm{tr}'\, x_i \pmod{J^2}. \]
Then we have
\[ (\sum_{i_j=i} e_{j})^2 =1, \hspace{1em} \sum_{i_j=k} e_{j} =0 \]
for any $k \neq i$.
\end{lem}
\textit{Proof.}
For any $s \in \C$, consider a representation $\rho_{8}: F_n \rightarrow \mathrm{SL}(2,\C)$ defined by
\[ \rho_{8}(x_r) = \begin{cases}
                  A, \hspace{1em} & \text{if} \hspace{1em} \hspace{0.5em} r = i, \\
                  E_2, & \text{if} \hspace{1em} \hspace{0.5em} r \neq i.
                 \end{cases}\]
Then from $\mathrm{tr}'(w) \equiv \mathrm{tr}'\, x_i \pmod{J^2}$ and Lemma {\rmfamily \ref{L-4}}, we obtain
\[ m^2 s + ({\text{terms of degree $\geq 2$}}) = s + ({\text{terms of degree $\geq 2$}}) \]
where $m = \sum_{i_j=i} e_{j}$. Hence $m^2=1$.

\vspace{0.5em}

Similarly, for any $k \neq i$ and $s \in \C$, considering a representation $\rho_{9}: F_n \rightarrow \mathrm{SL}(2,\C)$ defined by
\[ \rho_{9}(x_r) = \begin{cases}
                  A, \hspace{1em} & \text{if} \hspace{1em} \hspace{0.5em} r = k, \\
                  E_2, & \text{if} \hspace{1em} \hspace{0.5em} r \neq k.
                 \end{cases}\]
we obtain $\sum_{i_j=k} e_{j} =0$. This completes the proof of Lemma {\rmfamily \ref{L-4}}. $\square$

\vspace{0.5em}

From this lemma, we see that for any $\sigma \in \mathcal{E}_{F_n}(1)$ and $1 \leq i \leq n$,
\[ x_i^{\sigma} = x_i^{m_i} c_i, \hspace{1em} m_i= \pm 1 \]
for some $c_i \in \Gamma_{F_n}(2)$. Next, we show

\begin{lem}\label{L-5}
For any $\sigma \in \mathcal{E}_{F_n}(1)$, $m_1=m_2=\cdots=m_n$.
\end{lem}
\textit{Proof.}
Choose any $1 \leq i < j \leq n$. Consider a representation $\rho_{10}: F_n \rightarrow \mathrm{SL}(2,\C)$ defined by
\[ \rho_{10}(x_r) = \begin{cases}
                  A, \hspace{1em} & \text{if} \hspace{1em} \hspace{0.5em} r = i, j, \\
                  E_2, & \text{if} \hspace{1em} \hspace{0.5em} r \neq i, j.
                 \end{cases}\]
Then from $\mathrm{tr}'((x_i x_j)^{\sigma} ) \equiv \mathrm{tr}'\, x_i x_j \pmod{J^2}$, we see
\[ (m_i + m_j)^2 s + ({\text{terms of degree $\geq 2$}}) = 4s + ({\text{terms of degree $\geq 2$}}). \]
Hence we obtain $m_i = m_j$. This completes the proof of Lemma {\rmfamily \ref{L-5}}. $\square$

\vspace{0.5em}

Therefore we see that for any $\sigma \in \mathcal{E}_{F_n}(1)$ and $1 \leq i \leq n$,
\[ x_i^{\sigma} = x_i^{m_{\sigma}} c_i, \hspace{1em} m_{\sigma}= \pm 1 \]
for some $c_i \in \Gamma_{F_n}(2)$. Here assume
\[\begin{split}
   c_i & \equiv [x_2, x_1]^{e_{21}(i)} [x_3, x_1]^{e_{31}(i)} \cdots [x_{n}, x_{n-1}]^{e_{n \, n-1}(i)} \\
       & \hspace{3em} \cdot [x_2, x_1, x_1]^{e_{211}(i)} \cdots [x_n, x_{n-1}, x_n]^{e_{n \, n-1 \, n}(i)} \pmod{\Gamma_{F_n}(4)}
  \end{split}\]
for $e_{ba}(i), e_{bac}(i) \in \Z$. Here in the right hand side of the equation, terms $[x_b, x_a]$ for $b>a$ are multiplicated according to
the lexicographic ordering
\[ (b, a) < (b', a') \,\, \Longleftrightarrow \,\, \text{$a < a'$ or, $a=a'$ and $b < b'$}, \]
and terms $[x_b, x_a, x_c]$ for $b>a \leq c$ are multiplicated according to the lexicographic ordering
\[ (b, a, c) < (b', a', c') \,\, \Longleftrightarrow \,\, \begin{cases}
         a < a', \\
         a=a' \,\,\, \text{and} \,\,\, b < b' \,\,\, \text{or}, \\
         a=a', b=b' \,\,\, \text{and} \,\,\, c<c'.
       \end{cases}\]
We remark that
\[ \{ [x_b, x_a] \,|\, 1 \leq a < b \leq n \} \,\,\, \text{and} \,\,\, \{ [x_b, x_a, x_c] \,|\, b > a \leq c \} \] 
form basis of $\Gamma_{F_n}(2)/\Gamma_{F_n}(3)$ and $\Gamma_{F_n}(3)/\Gamma_{F_n}(4)$ as free abelian groups respectively.
(For details, see \cite{Hal}, for example.)

\begin{lem}\label{L-6}
As the notation above, $e_{ba}(i) = 0$ if $a, b \neq i$.
\end{lem}
\textit{Proof.}
We prove this lemma for the case where $a< b<i$. The other cases are proved in a similar way.
For any $s, t, u \in \C$, consider a representation $\rho_{11} : F_n \rightarrow \mathrm{SL}(2,\C)$ defined by
\[\rho_{11}(x_a) := \begin{pmatrix} 1 & s \\ 0 & 1 \end{pmatrix}, \,\,\, \rho_{11}(x_b) := \begin{pmatrix} 1 & 0 \\ t & 1 \end{pmatrix}, \,\,\,
    \rho_{11}(x_i) := \begin{pmatrix} 1-u & u \\ -u & 1+u \end{pmatrix}  \]
and $\rho_{11}(x_r)= E_2$ for $r \neq a, b, i$. Then for any $e \in \Z$, we have
{\small
\[\begin{split}
  \rho_{11}([x_b, x_a]^e) & \equiv \begin{pmatrix} 1-est & es^2t \\ -est^2 & 1+est \end{pmatrix}, \\
  \rho_{11}([x_i, x_a]^e) & \equiv \begin{pmatrix} 1+ e(su-su^2) & -e(2su+s^2u-su^2) \\ -e su^2 & 1-e(su-su^2) \end{pmatrix}, \\
  \rho_{11}([x_i, x_b]^e) & \equiv \begin{pmatrix} 1+e(tu+tu^2) & -etu^2 \\ e(2tu+tu^2+t^2u) & 1-e(tu+tu^2) \end{pmatrix},
  \end{split}\]
\[\begin{split}
  \rho_{11}([x_b, x_a, x_a]^e) & \equiv \begin{pmatrix} 1 & -2e s^2t \\ 0 & 1 \end{pmatrix}, \hspace{7em}
  \rho_{11}([x_b, x_a, x_b]^e) \equiv \begin{pmatrix} 1 & 0 \\ 2est^2 & 1 \end{pmatrix}, \\
  \rho_{11}([x_b, x_a, x_i]^e) & \equiv \begin{pmatrix} 1 & -e (s^2t+2stu) \\ -2e stu & 1 \end{pmatrix}, \hspace{1em}
  \rho_{11}([x_i, x_a, x_a]^e) \equiv \begin{pmatrix} 1 & 2e s^2u \\ 0 & 1 \end{pmatrix}, \\
  \rho_{11}([x_i, x_a, x_b]^e) & \equiv \begin{pmatrix} 1-2e stu & e (s^2u - su^2) \\ -2 e stu & 1 + 2 e stu \end{pmatrix}, \\
  \rho_{11}([x_i, x_a, x_i]^e) & \equiv \begin{pmatrix} 1+2e su^2 & -2e su^2 \\ 2e su^2 & 1-2e su^2 \end{pmatrix}, \hspace{2.3em}
  \rho_{11}([x_i, x_b, x_b]^e) \equiv \begin{pmatrix} 1 & 0 \\ -2e t^2u & 1 \end{pmatrix}, \\
  \rho_{11}([x_b, x_a, x_b]^e) & \equiv \begin{pmatrix} 1-etu^2 & 2e tu^2 \\ 0 & 1+e tu^2 \end{pmatrix}
\end{split}\]
}
modulo
\[ \{ F=(f_{ij}) \,|\, {\text{each of $f_{ij}$ is a polynomial of $s, t, u$ of degree $\geq 4$}} \}. \]
Hence from $\mathrm{tr}'(x_i)^{\sigma} \equiv \mathrm{tr}' x_i \pmod{J^2}$, we have
\[\begin{split}
   \mathrm{tr}' & \Big{(} \begin{pmatrix} 1-m_{\sigma} u & m_{\sigma} u \\ - m_{\sigma}u & 1+m_{\sigma}u \end{pmatrix}
    \rho_{11}([x_2, x_1]^{e_{21}(i)} [x_3, x_1]^{e_{31}(i)} \cdots [x_n, x_{n-1}, x_n]^{e_{n \, n-1 \, n}(i)}) \Big{)} \\
   & \hspace{1em} = 0 + ({\text{terms of degree $\geq 4$}}).
  \end{split}\]
Therefore, by observing the coefficient of $stu$, we obtain $2 m_{\sigma} e_{ba}(i) =0$, and hence $e_{ba}(i)=0$.
This completes the proof of Lemma {\rmfamily \ref{L-6}}. $\square$

\vspace{0.5em}

From this lemma, we see that for any $\sigma \in \mathcal{E}_{F_n}(1)$ and $1 \leq i \leq n$,
\[\begin{split}
   c_i & \equiv [x_i, x_1]^{e_{i1}(i)} [x_i, x_2]^{e_{i2}(i)} \cdots [x_{i}, x_{i-1}]^{e_{i \, i-1}(i)} \\
       & \hspace{3em} \cdot [x_{i+1}, x_i]^{e_{i+1 \, i}(i)} \cdots [x_n, x_i]^{e_{n \, i}(i)} \pmod{\Gamma_{F_n}(3)}.
  \end{split}\]
Next, we show
\begin{lem}\label{L-7}
As the notation above, for any $1 \leq i \leq n$,
\[ e_{i1}(1) = e_{i2}(2) = \cdots = e_{i \, i-1}(i-1) = -e_{i+1 \, i}(i+1) = \cdots = - e_{n i}(n). \]
\end{lem}
\textit{Proof.}
Choose any $a<b<i$, and fix it. We have
\[\begin{split}
   (x_a x_b)^{\sigma} & = x_a^{m_{\sigma}} c_a x_b^{m_{\sigma}} c_b \equiv x_a^{m_{\sigma}} x_b^{m_{\sigma}} c_a c_b \pmod{\Gamma_{F_n}(3)}, \\
    & \equiv x_a^{m_{\sigma}} x_b^{m_{\sigma}} [x_2, x_1]^{e_{21}(a) + e_{21}(b)} \cdots [x_n, x_{n-1}]^{e_{n \, n-1}(a) + e_{n \, n-1}(b)}
      \pmod{\Gamma_{F_n}(3)}.
  \end{split}\]
By the same argument as that in Lemma {\rmfamily \ref{L-6}, from an equation
$\mathrm{tr}' \rho_{11}((x_a x_b)^{\sigma}) \equiv \mathrm{tr}' \rho_{11}(x_a x_b) \pmod{J^2}$, we see
\[ st + 2 m_{\sigma}(-e_{ia}(a) - e_{ia}(b) + e_{ib}(a) + e_{ib}(b))stu \equiv st + ({\text{terms of degree $\geq 4$}}). \]
Hence we obtain $e_{ib}(b) = e_{ia}(a)$.

\vspace{0.5em}

Next, choose any $a<i<b$, and fix it. Similarly, from an equation
$\mathrm{tr}' \rho_{11}((x_a x_b)^{\sigma}) \equiv \mathrm{tr}' \rho_{11}(x_a x_b) \pmod{J^2}$,
we see
\[ st + 2 m_{\sigma}(-e_{ia}(a) - e_{ia}(b) - e_{bi}(a) - e_{bi}(b))stu \equiv st + ({\text{terms of degree $\geq 4$}}). \]
Hence we obtain $e_{bi}(b) = - e_{ia}(a)$.
This completes the proof of Lemma {\rmfamily \ref{L-7}}. $\square$

\vspace{0.5em}

For any $1 \leq i \leq n$, set $e_i :=  e_{i1}(1) = \cdots = - e_{n i}(n)$.
By Lemma {\rmfamily \ref{L-7}}, we see
\[ x_i^{\sigma} \equiv x_i^{m_{\sigma}} [x_1, x_i]^{e_1} \cdots [x_{i-1}, x_i]^{e_{i-1}} [x_{i+1}, x_i]^{e_{i+1}}
      \cdots [x_n, x_i]^{e_n} \pmod{\Gamma_{F_n}(3)}. \]
Then we show
\begin{lem}\label{L-8}
As the notation above, for any $\sigma \in \mathcal{E}_{F_n}(1)$, $m_{\sigma}=1$.
\end{lem}
\textit{Proof.}
Assume $m_{\sigma}=-1$. For any $1 \leq j \leq n$, let $\iota_j \in \mathrm{Inn}\,F_n$ be an inner automorphism of $F_n$ defined by
$x \mapsto x_j x x_j^{-1}$ for any $x \in F_n$. Then for any $1 \leq i \leq n$ and $e \in \Z$,
we have $x_i^{\iota_j^e} \equiv [x_j, x_i]^e x_i \pmod{\Gamma_{F_n}(3)}$. 
An element $\sigma' := \sigma \iota_1^{e_1} \cdots \iota_n^{e_n} \in \mathcal{E}_{F_n}(1)$ satisfies
\[ x_i^{\sigma'} \equiv x_i^{-1} \pmod{\Gamma_{F_n}(3)} \]
for each $1 \leq i \leq n$.
Hence, for any $a < b< i$, we see $(x_a x_b x_i)^{\sigma'} \equiv x_a^{-1} x_b^{-1} x_i^{-1} \pmod{\Gamma_{F_n}(3)}$.

\vspace{0.5em}

By an argument similar to that in Lemma {\rmfamily \ref{L-7}}, from an equation
$\mathrm{tr}' \rho_{11}((x_a x_b x_i)^{\sigma'}) \equiv \mathrm{tr}' \rho_{11}(x_a x_b x_i) \pmod{J^2}$,
we obtain
\[ st + tu -su +stu = st + tu - su - stu + ({\text{terms of degree $\geq 4$}}), \]
and hence the contradiction. This completes the proof of Lemma {\rmfamily \ref{L-8}}. $\square$

\vspace{0.5em}

As a corollary, we see
\begin{cor}
For any $n \geq 3$, $\mathcal{E}_{F_n}(1) \subset \mathrm{IA}_n$.
\end{cor}

\vspace{0.5em}

Now, we have
\begin{lem}\label{L-11}
For any $n \geq 3$, $\mathcal{A}_{F_n}(2) \subset \mathcal{E}_{F_n}(1)$
\end{lem}
\textit{Proof.}
For any $\sigma \in \mathcal{A}_{F_n}(2)$, and any $x \in F_n$, we have $x^{\sigma} = xy$ for some $y \in \Gamma_{F_n}(3)$.
Since $\Gamma_{F_n}(3)$ is generated by elemnts type of $[a,b,c]$ for $a, b, c \in F_n$, we can write
\[ y = [a_1, b_1, c_1]^{e_1} \cdots [a_r, b_r, c_r]^{e_r}, \hspace{1em} e_j= \pm 1. \]
Hence, using Lemma {\rmfamily \ref{L-10}} recursively,
we obtain $\mathrm{tr}' x^{\sigma} \equiv \mathrm{tr}' x \pmod{J^2}$ for any $x \in F_n$.
This completes the proof of Lemma {\rmfamily \ref{L-11}}.
$\square$

\vspace{0.5em}

Then we have
\begin{thm}\label{T-1}
For any $n \geq 3$, $\mathcal{E}_{F_n}(1) = \mathrm{Inn}\,F_n \cdot \mathcal{A}_{F_n}(2)$.
\end{thm}
\textit{Proof.}
Recall the argument in the former part of Lemma {\rmfamily \ref{L-8}}. For any $\sigma \in \mathcal{E}_{F_n}(1)$, there exists some
$\iota \in \mathrm{Inn}\,F_n$ such that
\[ x^{\sigma \iota} \equiv x \pmod{\Gamma_{F_n}(3)}  \]
for any $x \in F_n$. This shows that $\sigma \iota \in \mathcal{A}_{F_n}(2)$.
This completes the proof of Theorem {\rmfamily \ref{T-1}}. $\square$

\vspace{0.5em}

At the end of this subsection, we prove
\begin{thm}\label{T-2}
For any $k \geq 1$, $\mathcal{A}_{F_n}(2k) \subset \mathcal{E}_{F_n}(k)$.
\end{thm}
\textit{Proof.}
For any $\sigma \in \mathcal{A}_{F_n}(2k)$ and $x \in F_n$, we have $x^{\sigma} = x c$ for some $c \in \Gamma_{F_n}(2k+1)$.
By Lemma {\rmfamily \ref{l-1}}, the element $c$ is written as
\[ c = c_1^{e_1} c_2^{e_2} \cdots c_r^{e_r} \]
for some left-normed commutatiors $c_i$ of weight $2k+1$ and $e_i = \pm1$.
Hence, from Proposition {\rmfamily \ref{P-7}} and Corollary {\rmfamily \ref{C-P-7}},
we obtain
\[ \mathrm{tr}' x^{\sigma} \equiv \mathrm{tr}' x \pmod{J^{k+1}}. \]
This shows that $\sigma \in \mathcal{E}_{F_n}(k)$.
This completes the proof of Theorem {\rmfamily \ref{T-2}}. $\square$

\vspace{0.5em}

\subsection{Graded quotients $\mathrm{gr}^k(\mathcal{E}_{F_n})$}
\hspace*{\fill}\ 

\vspace{0.5em}

In this subsection, we study the graded quotients $\mathrm{gr}^k(\mathcal{E}_{F_n}) := \mathcal{E}_{F_n}(k)/\mathcal{E}_{F_n}(k+1)$.
Since each $\mathcal{E}_{F_n}(k)$ is a normal subgroup of $\mathrm{Aut}\, F_n$, the group $\mathrm{Aut}\, F_n$ naturally acts on 
$\mathrm{gr}^k(\mathcal{E}_{F_n})$ by the conjugation from the right. Furthermore since $\{ \mathcal{E}_{F_n}(k) \}$ is a central
filtration, the action of $\mathcal{E}_{F_n}(1)$ on $\mathrm{gr}^k(\mathcal{E}_{F_n})$ is trivial. Hence we can consider
each $\mathrm{gr}^k(\mathcal{E}_{F_n})$ as an $\mathrm{Aut}\,F_n/\mathcal{E}_{F_n}(1)$-module.
Here, we introduce Johnson homomorphism like homomorphisms to study the $\mathrm{Aut}\,F_n/\mathcal{E}_{F_n}(1)$-module structure of
$\mathrm{gr}^k(\mathcal{E}_{F_n})$.

\vspace{0.5em}

To begin with, for any $k \geq 1$ and $\sigma \in \mathcal{E}_{F_n}(k)$, define a map
$\eta_k(\sigma) : \mathrm{gr}^1(J) \rightarrow \mathrm{gr}^{k+1}(J)$ by
\[ \eta_k(\sigma)(f) := s_{\sigma}(f) = f^{\sigma} - f \in \mathrm{gr}^{k+1}(J) \]
for any $f \in J$. The well-definedness of the map $\eta_k(\sigma)$ follows from Lemma {\rmfamily \ref{L-2}}.
It is easily seen that $\eta_k(\sigma)$ is a homomorphism between abelian groups.

\vspace{0.5em}

Then we have a map $\eta_k : \mathrm{gr}^k(\mathcal{E}_{F_n}) \rightarrow \mathrm{Hom}_{\Z}(\mathrm{gr}^1(J), \mathrm{gr}^{k+1}(J))$
defined by $\sigma \mapsto \eta_k(\sigma)$. For any $\sigma, \tau \in \mathcal{E}_{F_n}(k)$, from (1) of Lemma {\rmfamily \ref{L-1}},
and from Lemma {\rmfamily \ref{L-2}}, we see
\[ s_{\sigma \tau}(f) =(s_{\sigma}(f))^{\tau} + s_{\tau}(f) \equiv s_{\sigma}(f) + s_{\tau}(f) \pmod{J^{k+2}}. \]
This shows that $\eta_k$ is a homomorphism of abelian groups. By the definition, each of $\eta_k$ is injective.
Furthermore, we have
\begin{lem}\label{L-12}
For each $k \geq 1$, $\eta_k$ is an $\mathrm{Aut}\,F_n/\mathcal{E}_{F_n}(1)$-equivariant.
\end{lem}
\textit{Proof.}
It suffices to show that $\eta_k$ is an $\mathrm{Aut}\,F_n$-equivariant.
For any $\sigma \in \mathrm{Aut}\,F_n$ and $\tau \in \mathcal{E}_{F_n}(k)$, we see
\[\begin{split}
   \eta_k(\tau \cdot \sigma)(f) & = \eta_k(\sigma^{-1} \tau \sigma)(f) = s_{\sigma^{-1} \tau \sigma}(f), \\
   (\eta_k(\tau) \cdot \sigma)(f) & = ( \eta_k(\tau)(f^{\sigma^{-1}}) )^{\sigma}
     = s_{\tau}(f^{\sigma^{-1}})^{\sigma} = ( f^{\sigma^{-1}\tau} - f^{\sigma^{-1}})^{\sigma} \\
               & = f^{\sigma^{-1}\tau \sigma} - f = s_{\sigma^{-1} \tau \sigma}(f)
  \end{split}\]
for any $f \in J$. Hence we have $\eta_k(\tau \cdot \sigma) = \eta_k(\tau) \cdot \sigma$.
This means $\eta_k$ is an $\mathrm{Aut}\,F_n$-equivariant homomorphism. 
This completes the proof of Lemma {\rmfamily \ref{L-12}}. $\square$

\vspace{0.5em}

Using the homomorphisms $\eta_k$, we see that $\mathrm{gr}^k(\mathcal{E}_{F_n})$ is an $\mathrm{Aut}\,F_n/\mathcal{E}_{F_n}(1)$-submodule
of the $\Q$-vector space $\mathrm{Hom}_{\Z}(\mathrm{gr}^1(J), \mathrm{gr}^{k+1}(J))$, and hence we obtain
\begin{thm}\label{P-4}
For any $n \geq 3$,
\begin{enumerate}
\item Each of $\mathrm{gr}^k(\mathcal{E}_{F_n})$ is torsion-free.
\item $\mathrm{dim}_{\Q} (\mathrm{gr}^k(\mathcal{E}_{F_n}) \otimes_{\Z} \Q) < \infty$.
\end{enumerate}
\end{thm}
As a corollary to Theorem {\rmfamily \ref{T-2}}, we see that $\mathrm{gr}^1(\mathcal{E}_{F_n})$ is finitely generated.
In general, however, it seems to be quite a difficult to determine the $\mathrm{Aut}\,F_n/\mathcal{E}_{F_n}(1)$-module
structure of $\mathrm{gr}^k(\mathcal{E}_{F_n})$ even the case where $k=1$.

\section{Acknowledgments}\label{S-Ack}

A part of this work was done when the authors stayed at Hokkaido University in the summer of 2011.
They would like to thank Professor Toshiyuki Akita for inviting us to Hokkaido University.

Both authors are supported by Grant-in-Aid for Young Scientists (B) by JSPS.

\end{document}